\documentclass[11pt,leqno,twoside]{amsart}
\usepackage{mathpazo}
\usepackage{geometry}
\usepackage{tikz-cd}
\usepackage{graphicx} 
\usepackage{amsmath}
\usepackage{amssymb}
\usepackage{extarrows}
\usepackage{csquotes}
\usepackage{comment}
\usepackage{mathrsfs}
\usepackage{amsthm}
\usepackage{enumitem}
\usepackage{color}
\usepackage{tikz}
\usepackage[colorlinks=true, linkcolor=blue, urlcolor=blue, citecolor=blue]{hyperref}
\newtheoremstyle{plainsl}
{8pt plus 2pt minus 4pt}
{8pt plus 2pt minus 4pt}
{\slshape}
{0pt}
{\bfseries}
{.}
{5pt plus 1pt minus 1pt}
{}
\theoremstyle{plainsl}
\newtheorem{theorem}{Theorem}[section]

\newtheorem{prop*}{Proposition}

\newtheorem{lemma}[theorem]{Lemma}
\newtheorem{corollary}[theorem]{Corollary}
\newtheorem{conjecture}[theorem]{Conjecture}
\newtheorem{remark}[theorem]{Remark}
\theoremstyle{definition}
\newtheorem{definition}[theorem]{Definition}

\def\ba#1\ea{\begin{align*}#1\end{align*}}

\newcommand{\C}{\mathbb{C}}

\newcommand{\Q}{\mathbb{Q}}
\newcommand{\R}{\mathbb{R}}

\newcommand{\cD}{\mathcal{D}}

\newcommand{\cH}{\mathcal{H}}

\newcommand{\cM}{\mathcal{M}}

\newcommand{\bD}{\mathbf{D}}

\newcommand{\sD}{\mathscr{D}}

\newcommand{\gr}{\mathrm{Gr}}

\newcommand{\DR}{\mathrm{DR}}

\newcommand{\MHM}{\mathrm{MHM}}

\newcommand{\mhm}{\mathrm{MHM}}
\newcommand{\td}{\mathrm{td}}

\title[Hodge index theorem for symmetric products]{\(S_n\)-equivariant \(L\)-classes of Cartesian powers of Witt spaces and the Hodge index conjecture for symmetric products}
\author{Mohammadali Aligholi}
\address{M. Aligholi: Department of Mathematics, University of Wisconsin-Madison, 480 Lincoln Drive, Madison WI 53706-1388, USA}
\email{aligholi@wisc.edu}
\date{}
\begin{document}
\raggedbottom
\begin{abstract}

For a compact oriented even-dimensional Witt pseudomanifold $X$, we compute the $S_n$-equivariant Banagl--Zagier $L$-classes of $X^n$ in terms of the Goresky--MacPherson $L$-classes of $X$ and the intersection cohomology Euler characteristic of $X$. This generalizes a theorem of \mbox{Zagier}  for \mbox{manifolds}. Our approach is based on analyzing the equivariant structure of the signature operator of $X^n$ using equivariant $KK$-theory and power operations in analytic $K$-homology. We establish an analytic $K$-homology analogue of Moonen's algebraic localization theorem by proving that the localized Chern character of the $n$-th power operation, evaluated at an $n$-cycle, agrees with the $n$-th homological Adams operation applied to the Chern character. As an application, we prove that the characteristic class Hodge index conjecture for singular varieties, formulated by Brasselet--Sch\"urmann--Yokura, is preserved under taking symmetric products of complex projective varieties.

\end{abstract}
\setcounter{tocdepth}{3}
\maketitle
\tableofcontents

\section{Introduction}

Let $G$ be a finite group and $M$ be an oriented
closed even-dimensional manifold equipped with an
orientation-preserving $G$-action. Let $g\in G$. The $g$-signature is the difference of the traces of $g$ on the positive and negative subspaces of the Hermitian middle-dimensional intersection pairing on $H^{\dim M/2}(M,\mathbb C)$.  Assume the components of the fixed point set $M^g$ are orientable. The
Atiyah--Singer $G$-signature theorem \cite{AS3} gives cohomological local Atiyah--Singer $L$-classes on the fixed point set  
\[L_{an}^*(M,g)\in H^*(M^g,\mathbb{C}) \]
whose evaluation on the fundamental class $[M^g]$ is the \(g\)-signature of $M$. We work with homological $L$-classes in this paper, which in this case $$L^{an}_*(M,g)=L_{an}^*(M,g)[M^g] \in H_*(M^g,\mathbb{C})$$ is the Poincaré dual to $L_{an}^*(M,g)$. The pushforward of $L^{an}_*(M,g)$ by the inclusion $j_g:M^g\hookrightarrow M$ is $L^{an}_*(g,M)\in H_*(M,\mathbb{C})$; the (homological) equivariant (Atiyah-Singer) $L$-classes of $M$ with respect to $g$, which was defined topologically by Zagier \cite{Z}.

Assume  $M$ is connected for simplicity.  The symmetric group $S_n$ acts on the product $M^n$ by permutation of coordinates. Let $\sigma_n\in S_n$ be an $n$-cycle. Then the fixed point set $(M^n)^{\sigma_n}$  is the diagonal and can be identified with  $(M^n)^{\sigma_n}\cong \Delta(M)\cong M$.  Zagier showed 
\cite[Prop.~2]{Z2}:
\begin{equation}\label{zagiersm}
    L^{an}_*(M^n,\sigma_n)=
\begin{cases}
\chi(M), & \text{if } n \text{ is even},\\[4pt]
\Psi_nL^{an}_*(M), & \text{if } n \text{ is odd}
\end{cases}
\end{equation}
where $\chi(M)$ is the Euler characteristic of $M$ viewed as an element of $H_0(M,\C)\cong \C$, and $L^{an}_*(M)$ is the (Atiyah-Singer) $L$-class of $M$, and $\Psi_n$ is the $n$-th homological Adams operation which acts on $H_{2k}$ by multiplication by $n^{-k}$.

If $X$ is a triangulated Witt pseudomanifold equipped with an orientation-preserving simplicial $G$-action, M. Banagl \cite{B} generalized the result of Zagier and  constructed equivariant homological $L$-classes $$L_*(g,X)\in H_*(X,\mathbb{C})$$ using topological methods analogous to those in \cite{Z}, which recover the Goresky--MacPherson $L$-classes of $X$ when $g=1$. 

The author proved in \cite{A} that one can obtain $L_*(g,X)$ by pushforward of some localized classes $L^{an}_*(X,g)\in H_*(X^g,\mathbb{C})$, obtained as the localized Chern character of the equivariant analytic $K$-homology class of the signature operator $\left[D_X^{sign,G}\right]$ defined by Banagl--Leichtnam--Piazza \cite{BLP}. More precisely, if $j_g:X^g\hookrightarrow X$ is the inclusion, then
\begin{equation} \label{realtionbanaglme}
    L_*(g,X)=j_{g*}(\Psi_2(L^{an}_*(X,g)))
\end{equation}
\begin{equation}\label{loclaizedlclasschern}
    L^{an}_*(X,g):=Ch(\left[D_X^{sign,G}\right],g)
\end{equation}
 where  $$Ch(-,g):K^{G,an}_0(X)\rightarrow H_*(X^g,\mathbb{C})$$ is the localized Chern character  transformation (Definition \ref{localizedchernclass} and \cite{A}). For $g=1$, $ j_{g*}L^{an}_*(X,1)=L^{an}_*(X)$ is the \textit{unstable} Atiyah–Singer class,  where $\Psi_2L^{an}_*(X)=L_*(X)$ is the Goresky--MacPherson $L$-class of $X$.

 We generalize (\ref{zagiersm}) to Witt pseudomanifolds:
\begin{theorem}\label{equianvtLclasszagiersingular}
    Let $X$ be a compact oriented connected even-dimensional smoothly stratified Witt pseudomanifold.   Let $\sigma_n\in S_n$ be an $n$-cycle. Then under the identification $(X^n)^{\sigma_n}\cong \Delta(X)\cong X$ 
    \[
L^{an}_*(X^n,\sigma_n)=
\begin{cases}
I\chi(X), & \text{if } n \text{ is even},\\[4pt]
\Psi_nL^{an}_*(X), & \text{if } n \text{ is odd}
\end{cases}
\]
where $I\chi(X)$ is the Euler characteristic of the intersection cohomology of $X$ viewed as an element of $H_0(X,\C)\cong \C$.
\end{theorem}

By the multiplicativity of localized Chern character of signature operator  (see Theorem \ref{symmetricfactorchern}), for an arbitrary $\sigma\in S_n$ we get
\begin{equation}
    L_*^{an}(X^n,\sigma)=\prod_{r,j} L^{an}_*(X^r,\sigma_{r,j}) \in  H_*(X^{\sum_r k_r},\C) \cong H_*((X^n)^{\sigma},\C)
\end{equation}
 where $\sigma=\prod_r\prod^{k_r}_{j=1}\sigma_{r,j}$ is the cycle decomposition and $\sigma_{r,j}$ is an $r$-cycle, and $\prod_{r,j}$ denotes the external
product of homology classes.

Zagier proved (\ref{zagiersm}) using the Atiyah–Singer fixed-point formula \cite{AS3}, which expresses the (cohomological) localized class in terms of a characteristic class factor of the normal bundle of the fixed-point set. An analogous localization formula was obtained in \cite{BLP} in the Witt case under a normally nonsingular condition on the inclusion of fixed-point sets. However, in the case of the action of $S_n$ on $X^n$, the inclusion of the diagonal is not normally nonsingular in general, so we cannot use this method to prove Theorem \ref{equianvtLclasszagiersingular}.
The approach taken in this paper is different and is based on
 studying the equivariant structure of the signature operator of $X^n$ and power operations in $K$-homology. In the course of the proof, we establish several results that are also of independent interest.

To continue the discussion we need to give some background on analytic $K$-homology in the sense of Kasparov \cite{K}.
Let $X$ be a compact oriented even dimensional smoothly stratified Witt pseudomanifold, see \cite{ALMP, ABP1} for terminology, e.g., a smoothly stratified space has the same data as a Thom-Mather stratified space. An analytic $K$-homology class of the signature operator $$\left[D_X^{sign}\right]\in K_{0}^{an}(X):=KK_0(C(X),\C)$$ was constructed by Moscovici--Wu \cite{MW}, based on the work of Cheeger \cite{Ch}, and later by P. Albin, E. Leichtnam, R. Mazzeo, and P. Piazza \cite{ALMP}. An even unbounded Kasparov $K$-cycle defining a class in $K_{0}^{an}(X)$, in the sense of \cite{BJ}, consists of a tuple $$(\mathcal H,\Gamma,\rho,\bD)$$ where  $\cH$ is  a $\mathbb{Z}_2$-$\Gamma$-graded Hilbert space with a \(*\)-representation $\rho$ of the algebra of continuous functions of $X$ and a self-adjoint densely defined odd operator $\bD$ satisfying some properties.   The signature operator $\bD=D_X^{sign}$  has a domain inside  $$\mathcal{H}=L^2(X_{reg},\Omega^*):=L^2(X_{reg},\Omega^*,g)$$ i.e., the Hilbert space of $L^2$ sections of differential forms on the regular part
 endowed with a metric $g$ quasi-isometric to a piecewise-flat conic metric, such as a wedge metric considered in \cite{ALMP}. The grading on $\cH$ is the signature grading $\tau$ on differential forms defined by the Hodge $\ast$ operator, and the representation $\rho$ is the multiplication by functions. The resulting $K$-homology class is independent of the choice of metric. There is a homological Chern character $Ch:K_*^{an}(X)\rightarrow H_*(X,\Q)$ which by a result of  Moscovici--Wu \cite{MW}, 
\begin{equation}\label{mw_lclass}
    Ch(\left[D_X^{sign}\right])=L^{an}_*(X)
\end{equation}
where $\Psi_2L^{an}_*(X)$ is  the Goresky--MacPherson $L$-classes $L_*(X)$. Similarly, if $\cH$ is equipped with  the  parity-of-degree  grading $\epsilon$, the de Rham operator (Euler characteristic operator) $\bD=D^{dR}_X$ gives rise to a class \cite{A2} $$\left[D_X^{dR}\right]\in K_0^{an}(X)$$
 As we mentioned, the $K$-homology class of the signature operator is very rich, unfortunately, the $K$-homology class of the de Rham operator is as trivial as it can get, in the sense that it contains only one datum, its numerical index, which is the Euler characteristic. This is proved for manifolds by J. Rosenberg \cite{R2}, and generalized by the author to Thom–Mather stratified pseudomanifolds in \cite{A2}. We give the result here for Witt spaces. This is analogous to the algebraic result (\ref{mhc-1trivial}) mentioned below.

  \begin{theorem}\cite{A2}\label{maintheorem} 
    Let $X$ be a compact connected  oriented smoothly stratified Witt pseudomanifold. Then
    $[D^{dR}_X]\in K^{an}_0(X)$ is in the image of $ K^{an}_0(pt)\rightarrow K^{an}_0(X)$ induced by the inclusion of a  point of $X$. This implies
    \begin{equation}\label{cheuler}
Ch([D^{dR}_X])=Ch_0([D^{dR}_X])=index([D^{dR}_X])=I\chi(X) \in H_0(X,\mathbb{C})
    \end{equation}
    where   $I\chi(X)$ is the Euler characteristic of the intersection cohomology of $X$ viewed as an element of $H_0(X,\C)$.
\end{theorem}

Now consider the permutation action of $S_n$ on the coordinates of $X^n$.  If $X$ is an even dimensional Thom-Mather Witt pseudomanifold, then $X^n$ is an $S_n$-Thom-Mather Witt pseudomanifold, in the sense of \cite{BLP}, see Theorem \ref{Sn-mather-X^n}. Then by a result of loc. cit. there is an equivariant signature operator defined by considering an $S_n$-equivariant wedge metric, which in this case one can take the product metric.  Then this gives an equivariant $K$-homology class $\left[D^{sign,S_n}_{X^n}\right]\in K_0^{S_n,an}(X^n)$.

\begin{theorem}\label{ch_sign_x^n-zagier}  Let $X$ be an even-dimensional compact  oriented  smoothly stratified Witt pseudomanifold.   Then $X^n$ is an $S_n$-Thom-Mather space where  $S_n$ acts naturally by permuting the coordinates. Let $\sigma_n\in S_n$ be an $n$-cycle. Then
\[
Ch(\left[D^{sign,S_n}_{X^n}\right],\sigma_n)
=
\begin{cases}
\Psi_nCh([D^{dR}_X]), & \text{if } n \text{ is even},\\[4pt]
\Psi_nCh([D^{sign}_X]), & \text{if } n \text{ is odd}.
\end{cases}
\]
\end{theorem}
 Combined with (\ref{loclaizedlclasschern}), (\ref{mw_lclass}), and (\ref{cheuler}) this implies Theorem \ref{equianvtLclasszagiersingular}.

To obtain Theorem \ref{ch_sign_x^n-zagier}, we identify the $S_n$-equivariant structure of $\left[D^{sign,S_n}_{X^n}\right]$ by a careful analysis of gradings. Note that non-equivariantly, by a result of Albin-Banagl-Piazza \cite{ABP1}, we have
 \[[D^{sign}_{X^n}]=[D^{sign}_{X}]^{\boxtimes n}\]
  where $\boxtimes n$ is the $n$-fold external Kasparov product.

There is a natural $S_n$-action on the $n$-fold external Kasparov product, by graded permutation of elements. By this action, the $n$-fold external Kasparov product defines the power operation in $K$-homology (the analogue of Atiyah's power operation in $K$-theory, \cite{At}), see Section \ref{poweropsection}.
\begin{equation}\label{intropower}
    P_n:K^{an}_0(X)\rightarrow K_0^{S_n,an}(X^n)
\end{equation}
In general, 
$
    [D^{sign,S_n}_{X^n}] \not = P_n([D^{sign}_{X}])
$.
This boils down to the difference between the signature grading $\tau$ and the usual parity-of-degree grading $\epsilon$. For simplicity of notation let
 $ \left[D_X^{sign}\right]=[(\mathcal{H},
\tau,\rho,
D)]$.
If $\dim_{\R} X_{}=2m$, the signature grading  $\tau$ is defined on $k$-forms by \cite{AS3}
 \begin{equation}\label{siggradingdef}
     \tau|_{\Omega^k(X_{reg})}=i^{k(k-1)+m}\ast
 \end{equation}
 where $\ast:\Omega^k(X_{reg})\rightarrow \Omega^{2m-k}(X_{reg})$ is the Hodge $\ast$ operator and $\epsilon|_{\Omega^k(X_{reg})}=(-1)^k$. By basic calculation
 $\epsilon \tau=\tau\epsilon$. Define
 $\kappa:=\epsilon \tau$ where  $\kappa^2=1$.
 Let \[\mathcal{H}=\mathcal{H}_{0} \oplus \mathcal{H}_{1}, \quad \kappa|_{\mathcal{H}_{0}}=1, \kappa|_{\mathcal{H}_{1}}=-1\] in other words, on $\mathcal{H}_{0}$, $\tau, \epsilon$ coincide while on $\mathcal{H}_{1}$ they are opposite. We have $$D\epsilon=-\epsilon D,\quad D\tau=-\tau D,\quad \kappa D=D\kappa$$ thus $D$ preserves $Dom(D)\cap \mathcal{H}_{a}$ and defines for $a\in \{0,1\}$ (unbounded) Kasparov cycles
 \[S_a:=(\mathcal{H},
\tau,\rho,
D)|_{\mathcal{H}_{a}}:=(\cH_a,\tau_a,\rho_a,D_a)\]
where $\tau_a:=\tau|_{\cH_a}, \rho_a:=\rho|_{\cH_a}, D_a:=D|_{Dom(D)\cap \cH_a}$.  We have
\[[S_0]+[S_1]=\left[(\cH_0 \oplus \cH_1,\tau,\rho, D) \right]=[D^{sign}_{X}]\]
\begin{equation}
    [S_0]-[S_1]=\left[(\cH_0 \oplus \cH_1,\tau_0 \oplus (-\tau_1),\rho, D )\right]=\left[(\cH_0 \oplus \cH_1,\epsilon,\rho, D) \right]=[D^{dR}_{X}]
\end{equation}
To simplify the computations define $\widetilde S_0=S_0,~ {\widetilde S_1}=(\cH_1,-\tau_1,\rho_1,D_1)$, so $$[\tilde{S}_0]=[S_0],\quad[\widetilde S_1]=-[S_1],\quad  [D^{sign}_{X}]=[\widetilde S_0] -[\widetilde S_1],\quad [D^{dR}_{X}]=[\widetilde S_0] +[\widetilde S_1] $$
Note that  $\tau_a=(-1)^a\epsilon_a$.

\begin{theorem}\label{introtheoreminductform}
    Let $X$ be an even-dimensional oriented closed smoothly stratified Witt pseudomanifold. Then in $K^{S_n,an}_0(X^n)$ we have
    \begin{equation}\label{introindrel}
        \left[D^{sign,S_n}_{X^n}\right]=\sum^n_{k=0}(-1)^kInd_{S_{n-k}\times S_k}^{S_n}\left(P_{n-k}([\widetilde S_0])\boxtimes P_{k}([\widetilde S_1])\right)
    \end{equation}
    where $\left[D^{sign,S_n}_{X^n}\right]$ is the $S_n$-equivariant signature class in the sense of \cite{BLP}, $P_i$ is the $i$-th power operation in analytic $K$-homology (\ref{intropower}), and 
      $$Ind_{S_{n-k}\times S_k}^{S_n}:K_0^{S_{n-k}\times S_k,an}(X^{n})\rightarrow K_0^{S_{n},an}(X^{n})$$
    is the induction on the level of $KK$-theory mentioned in Section \ref{inductkcycles}.
\end{theorem}
\begin{remark}
    The motivation to find the relation  (\ref{introindrel}) is the following relation given by Atiyah \cite[proof of Prop. 2.3]{At} 
   \[ ([E]-[F])^{\otimes n}=
\sum_{k=0}^n (-1)^k
\operatorname{Ind}_{S_{n-k}\times S_k}^{S_n}
\left(
E^{\otimes(n-k)}\otimes F^{\otimes k}
\right)
\]
where $E,F$ are vector bundles and $\otimes^ n:K^0(X)\rightarrow K_{S_n}^0(X)\cong K^0(X)\otimes R(S_n)$, and $E^{\otimes(n-k)}$ is an $S_{n-k}$ vector bundle via the standard permutation while $S_k$ acts on $F^{\otimes k}$ via permutation and signs. 
\end{remark}

To compute the localized Chern character transformation $Ch(\left[D^{sign,S_n}_{X^n}\right],\sigma_n)$, we first prove that the localized Chern character of $n$-th power operation with respect to an
$n$-cycle is the same as the $n$-th homological Adams operation applied to the Chern character.

\begin{theorem} \label{analyticmoonen}
      Let $X$ be a compact CW-complex. $S_n$ acts on $X^n$ by permutation. Let $\sigma_n\in S_n$ be an $n$-cycle. Identify the fixed point set $(X^n)^{\sigma_n}=\Delta(X) \cong X$. For $\mathcal{F}\in K_0^{an}(X)$, the $n$-th power operation
defines a class   $ P_n(\mathcal{F})=\mathcal{F}^{\boxtimes n} \in K_0^{{S_n},an}(X^n)$.  Then under the identification $H_*((X^n)^{\sigma_n},\mathbb{C})\cong H_*(X,\mathbb{C})$, we have
\begin{equation}\label{eqanalyticmoonen}
    Ch(P_n(\mathcal{F}),\sigma_n)=\Psi_n Ch(\mathcal{F})
\end{equation} 
where $\Psi_n$ is the $n$-th homological Adams operation. 
\end{theorem}

Note that this is an analytic analogue of the following algebraic localization theorem, in the sense that the Chern character in analytic $K$-homology  corresponds to the Baum--Fulton-MacPherson Todd class transformation \cite{BFM} on the Grothendieck group of coherent sheaves $$\td_*:G_0(X):=K_0(Coh(X))\rightarrow H^{BM}_*(X,\mathbb{C})$$ The following was proved by B. Moonen in the case of $\mathcal{G}=\mathcal{O}_X$ using an embedding of a quasi-projective variety $X$ into a smooth complex algebraic variety $M$. It is mentioned in \cite{CMSSY} that this argument applies directly to the more general context of $G_0(X)$.
\begin{theorem} (\cite{M},  \cite[Lemma 3.6]{CMSSY})\label{algebraic_moonen} Let \(X\) be a quasi-projective variety, and  $\sigma_n$ be an $n$-cycle and $\mathcal{G}\in G_0(X)$. Then under the identification $H^{BM}_*((X^n)^{\sigma_n},\mathbb{C})\cong H^{BM}_*(X,\mathbb{C})$, we have
\begin{equation}\label{moonentdeq}
    \td_{*}(\mathcal{G}^{\boxtimes n},\sigma_n)=\Psi_n\td_{*}(\mathcal{G})
\end{equation}
where $\Psi_n$ is the $n$-th homological Adams operation. 
\end{theorem}

To prove Theorem \ref{analyticmoonen}, we first establish the result for a compact even-dimensional $spin^c$ manifold $M$ (Theorem \ref{problemmoon}), using a $spin^c$ version of the Atiyah--Singer localization theorem \cite{AS3} and Atiyah's power operations in $K$-theory \cite{At}. We then extend the result to analytic $K$-homology via the isomorphism between Baum--Douglas $K$-homology \cite{BD} and analytic $K$-homology.

Finally, we give an application of Theorem \ref{equianvtLclasszagiersingular} in complex  algebraic geometry. We show that if a conjectural characteristic  Hodge index theorem for singular varieties, explained below, holds for a complex projective $X$, then it holds for symmetric products $X^{n}/S_n$ as well.

Assume $X$ is a compact complex algebraic variety. Then Brasselet–Sch\"urmann–Yokura \cite{BSY} formulated a conjectural characteristic class analogue of the Hodge index theorem for singular complex algebraic varieties, which expresses the Goresky--MacPherson homology $L$-classes in terms of suitable Hodge-theoretic $L$-classes that we now explain.

Let $\MHM(X)$ denote the abelian category of algebraic mixed Hodge modules on $X$ \cite{Sa}. The \textit{unnormalized} Hodge-theoretic {\it Hirzebruch class transformation} \cite{BSY} is defined by the composition
\begin{align*}
MHT_{y\ast}: K_0(\mhm(X)) &\to H_*(X) \otimes \Q[y^{\pm 1}, (1+y)^{-1}] \\
MHT_{y\ast} &=\td_{\ast}  \circ MHC_y
\end{align*}
where $\td_{\ast}$
is the scalar extension of the
Baum--Fulton-MacPherson Todd class transformation $\td_*$, and a motivic Chern class transformation
\[
 MHC_y: K_0(\MHM(X)) \to  G_0(X)[y,y^{-1}]
\]
If $X$ is smooth and $M\in \MHM(X)$, let $(\cM,F_{\bullet}\cM)$ be the underlying filtered $\sD_X$-module of $M$, and set
\[
\begin{split}
MHC_y[M] &= \sum_{p}  \left[\gr^F_{-p}\DR(\cM)\right] \cdot(-y)^p  \\
&= \sum_{p,i} (-1)^i \left[\cH^i \gr^F_{-p}\DR(\cM)\right] \cdot(-y)^p \in G_0(X)[y,y^{-1}]
\end{split}
\]
with $\gr^F_{-p}\DR(\cM)$ the graded parts of the de Rham complex of $\cM$ with respect to the induced Hodge filtration. 

Let $IC'^H_X:=IC_X^H[-\dim_{\mathbb{C}}(X)]$ be the (shifted) intersection Hodge module and $$IT_{1*}(X):=MHT_{1*}([IC'^H_X]) \in H_*(X,\Q)$$
The following conjecture is attributed to Cappell and Shaneson.

\begin{conjecture} \label{BSYi}\cite{BSY}
Let $X$ be a compact pure-dimensional complex algebraic variety. Then
\begin{equation}\label{conji}
L_*(X)=\Psi_2 IT_{1*}(X)
\end{equation}
\end{conjecture}

Conjecture \ref{BSYi} holds in the smooth context by the work of Hirzebruch on the generalized Hirzebruch-Riemann-Roch theorem \cite{H}. The conjecture has been proved in many cases, including isolated singularities, rational homology manifolds, compact toric varieties, Schubert varieties, and complex algebraic surfaces and threefolds, see \cite{AMS} for an overview.

In \cite{A}, the author investigated the case of a finite quotient space $X/G$. Conjecture \ref{BSYi} holds for $X/G$ if we either assume that an equivariant $K$-theoretical version of the conjecture holds for $X$, or provided that it holds for all fixed point sets $X^g$ and that a suitable normally nonsingular inclusion assumption is satisfied. For a symmetric product $X^{(n)}=X^n/S_n$, the fixed point sets are products of diagonals. The diagonal map is not normally nonsingular in general; therefore, the argument given in loc. cit. does not apply.

\begin{theorem}\label{bsyconjforsymmetric}
     Let $X$ be a pure-dimensional complex projective variety. If Conjecture \ref{BSYi} holds for $X$, then it also holds for symmetric products $X^{(n)}=X^n/S_n$ for any $n\in \mathbb{N}$.
 \end{theorem}

 Theorem \ref{bsyconjforsymmetric} follows from an analogous argument given in \cite{CMSSY,CMSS}, by expressing characteristic classes of the quotient space $X^n/S_n$ as an average of localized equivariant classes on $X^n$ for both $MHT_{1*}$ and $L$-classes. The algebraic results for $MHT_{1*}([IC'^H_{X^{(n)}}])$ are proved in \cite{CMSSY} as follows. Note that 
 \[IC'^H_{X^{(n)}}=(IC'^H_X)^{(n)} \in D^b\mhm(X^{(n)})\]
where $$(IC'^H_X)^{(n)}:=\left(\pi_{n*}(IC'^H_X)^{\boxtimes n}\right)^{S_n}$$ is the $n$-th symmetric power of $IC'^H_{X}$ with an $S_n$ action defined by Maxim--Saito--Sch\"urmann \cite{MSS}.

 \begin{theorem}\cite[Lemma 3.1, 3.2, 3.3]{CMSSY} \label{it1symmetriccmssythm}
 Let $X$ be a pure-dimensional complex quasi-projective variety. Then
\begin{equation}\label{symfactorIT1class} 
    IT_{1*}(X^{(n)})=\frac{1}{n!}\sum_{\sigma\in S_n}  \prod_{r,j} (\pi^{\sigma_{r,j}}_*\Psi_rIT_{(-1)^{r+1}*}(X)) \in H^{BM}_*(X^{(n)},\mathbb{C})
\end{equation}
    where $\sigma=\prod_r\prod^{k_r}_{j=1}\sigma_{r,j}$ is the cycle decomposition, $\sigma_{r,j}$ is an $r$-cycle with $|\sigma_{r,j}|=r$ and
    \[\pi^{\sigma_{r,j}}:X\cong (X^{|\sigma_{r,j}|})^{\sigma_{r,j}} \hookrightarrow X^{|\sigma_{r,j}|} \rightarrow X^{(|\sigma_{r,j}|)}\]
with \[\prod_{r,j} X^{(|\sigma_{r,j}|)}\rightarrow X^{(n)}\]
induced by the Pontryagin product.
\end{theorem}
When $r$ is even, by \cite[Prop.~5.21]{S} and \cite{CMSSY} (assume $X$ is compact and connected)
\begin{equation}\label{mhc-1trivial}
    IT_{(-1)^{r+1}*}(X)=IT_{-1*}(X):=MHT_{-1*}([IC'^H_X])=c_0([IC'_X])=I\chi(X)
\end{equation}
is the degree-zero MacPherson Chern class of the underlying intersection cohomology complex, i.e., the corresponding Euler characteristic  under the identification
$H_0(X)\cong \Q$.
(Note that \cite{S} uses normalized Hirzebruch classes, whereas \cite{CMSSY} uses unnormalized Hirzebruch classes, as in the present paper.)
 Therefore, there are no higher characteristic classes besides the degree-zero part. This is the analogue of Theorem \ref{maintheorem}, the Chern character of the de Rham operator is only the degree-zero Euler characteristic of intersection cohomology. 

By a result of Banagl \cite[Thm. 10.1]{B}, one can express the $L$-classes of the quotient as an average of equivariant $L$-classes of $X^n$. Using Theorem \ref{equianvtLclasszagiersingular}, we can write the analogous averaging formula for $L$-classes:
\begin{theorem}\label{symLclassthm}
  Let $X$ be a  compact connected oriented even-dimensional smoothly stratified Witt pseudomanifold, then
\begin{equation}\label{symfactorLclass}
    L^{an}_*(X^{(n)})=\frac{1}{n!}\sum_{\sigma\in S_n}  \prod_{r,j} (\pi^{\sigma_{r,j}}_*\Psi_r C_{(-1)^r}(X)) \in H_*(X^{(n)},\mathbb{C})
\end{equation}
where $C_{-1}(X):=L_*^{an}(X)$, and $C_{1}(X):=I\chi(X)$.
\end{theorem}

\medskip
{\bf Acknowledgments.}  I thank my advisor, Lauren\c{t}iu Maxim, for his continued support and many helpful discussions.

\section{Background}\label{backgroundnotations}
Notation:

\begin{enumerate}
    \item For a finite group $G$, $R(G)$ denotes the complex representation ring of $G$. For an $R(G)$-module $M$, and $g\in G$, we denote by $M_g$ the localization of $M$ with respect to the prime ideal $\mathfrak{m}_g$ of $R(G)$ consisting of virtual representations whose character vanishes at $g$. Let
$\chi_g: R(G)_g \rightarrow \mathbb{C}$ defined by $\frac{\psi}{\zeta}\mapsto\frac{\psi(g)}{\zeta(g)}$.

\item $K^{top}_{0}, K^{an}_{0}, K_0^{geo}, G_0$  denote, respectively,  topological complex $K$-homology,  Kasparov $K$-homology $K^{an}_{0}(X)=KK_{0}(C_0(X),\mathbb{C})$ \cite{K},  Baum--Douglas $K$-homology \cite{BD}, and  $G_0(X)=K_0(Coh(X))$ is the Grothendieck group of algebraic coherent sheaves. We represent equivariant $KK$-theory by $KK^G(-,-)$.
\item The subscript $(n)$ in a homological transformation, such as $Ch_{(n)}$ or $\td_{n*}$, denotes composition with the $n$-th homological Adams operation $\Psi_n$ where  $\Psi_n|_{H_{2k}}=n^{-k}\operatorname{id}$. Similarly $\Psi^n$ denotes $n$-th cohomological Adams operation $\Psi^n|_{H^{2k}}=n^k\operatorname{id}$. Let $\psi^n$ denote $n$-th Adams operation in $K^0(M)$ \cite{At}.
\item The subscript $g$ denotes composition with the evaluation map $\chi_g:R(G)_g\rightarrow\mathbb{C}$. 

\end{enumerate}

\subsection{Equivariant unbounded Kasparov cycles}

The unbounded formalism for the Kasparov KK-group \cite{BJ} is as follows, see also \cite{Bl}. Let $X$ be a second-countable locally compact Hausdorff space.

\begin{definition}\label{kkdef}\cite{BJ}
An even unbounded $K$-cycle defining a class in
\[
K_0^{an}(X)=KK_0(C_0(X),\C)
\]
consists of a tuple $(\mathcal H,\Gamma,\rho_X,\bD)$, where
\begin{enumerate}[label=(\alph*)]
    \item $\mathcal H$ is a $\mathbb Z_2$-graded Hilbert
    space with grading operator $\Gamma$ satisfying
    $\Gamma^2=1$ and $\Gamma^*=\Gamma$.

    \item $\rho_X:C_0(X)\to\mathcal B(\mathcal H)$ is
    an even, nondegenerate $\ast$-representation, i.e.,
    \[
    \Gamma\rho_X(f)=\rho_X(f)\Gamma,
    \qquad f\in C_0(X)
    \]

    \item $\bD:\operatorname{Dom}(\bD)\subset\mathcal H
    \to\mathcal H$ is an odd self-adjoint unbounded
    operator with respect to $\Gamma$, i.e.,
    \[
    \Gamma\bD=-\bD\Gamma
    \]
\end{enumerate}
Furthermore, the following conditions are required:
\begin{enumerate}
    \item There exists a dense $\ast$-subalgebra
    $A\subset C_0(X)$ such that
    $\operatorname{Dom}(\bD)$ is invariant under
    $\rho_X(f)$ and $[\bD,\rho_X(f)]$ extends
    to a bounded operator for every $f\in A$.

    \item For every $f\in C_0(X)$,
    \[
    \rho_X(f)(1+\bD^2)^{-1}
    \]
    is a compact operator on $\mathcal H$.
\end{enumerate}
The bounded transform
\[
F=\bD(1+\bD^2)^{-1/2}
\]
defines an even \textit{bounded} Kasparov cycle in the sense of \cite{K} and hence a class
$[\bD]\in K_0^{an}(X)$.

When $X$ is compact, choosing $f=1\in C_0(X)$, the second condition implies and is equivalent to that $\bD$ has compact resolvent, which is equivalent to $\bD$ having purely discrete spectrum with finite multiplicities. 
\end{definition}

There is a notion of \textit{external Kasparov product}  \cite{K}
 \begin{equation}\label{external_kasparov_product}
     \boxtimes:K^{an}_0(X)\otimes K^{an}_0(Y)\rightarrow K^{an}_{0}(X\times Y)
 \end{equation}
which can be expressed more simply  in the unbounded model of $K$ homology as follows.  For two unbounded cycles 
\[[x]=[\mathcal{H}_X,\Gamma_X,\rho_X,\bD_X] \in K^{an}_0(X),\quad [y]=[\mathcal{H}_Y,\Gamma_Y,\rho_Y,\bD_Y] \in K^{an}_0(Y)\]
it is defined by
 \[[x]\boxtimes [y]=\left[
\mathcal{H}_X\widehat{\otimes} \mathcal{H}_Y,\Gamma_X\otimes \Gamma_Y,\,
\rho_X\widehat{\otimes}\rho_Y,\,
\bD_X\widehat{\otimes}1+1\widehat{\otimes}\bD_Y
\right]\]
here $\mathcal H_X\widehat\otimes\mathcal H_Y$ denotes
the graded Hilbert tensor product: the completion of
$\mathcal H_X\odot\mathcal H_Y$ for the inner product
determined by
\[
\langle \xi\otimes\eta,\xi'\otimes\eta'\rangle
=
\langle\xi,\xi'\rangle_{\mathcal H_X}
\langle\eta,\eta'\rangle_{\mathcal H_Y}
\]
equipped with the grading $\Gamma_X\otimes\Gamma_Y$
\cite[\S~14.4]{Bl}.

\[
\begin{aligned}
(\mathcal H_X\widehat{\otimes}\mathcal H_Y)^0
&=(\mathcal H_X^0\otimes\mathcal H_Y^0)
\oplus(\mathcal H_X^1\otimes\mathcal H_Y^1)\\
(\mathcal H_X\widehat{\otimes}\mathcal H_Y)^1
&=(\mathcal H_X^0\otimes\mathcal H_Y^1)
\oplus(\mathcal H_X^1\otimes\mathcal H_Y^0)
\end{aligned}
\]
where $(-)^{0,1}$ denotes the grading with respect to $\Gamma_{X(Y)}$,
and $\bD_X\widehat{\otimes}1 (1\widehat{\otimes}\bD_Y)$ are the graded tensor products of operators, i.e., if $A,B$ are homogeneous operators of degrees $|A|$ and $|B|$ (where an odd operator has degree \(1\), and an even operator has degree \(0\))
\[A\widehat{\otimes}B(\zeta\otimes \gamma)=(-1)^{|B||\zeta|} A\zeta\otimes B\gamma\]
therefore for odd operators $\bD_{X(Y)}$, under the identification $\mathcal H_X\widehat{\otimes}\mathcal H_Y=\mathcal H_X{\otimes}\mathcal H_Y$ we have
\[\,
\bD_X\widehat{\otimes}1+1\widehat{\otimes}\bD_Y={\bD_X\otimes 1 +\Gamma_X\otimes \bD_Y}\]
here the tensor-sum operator denotes its self-adjoint
closure from the initial domain
$\operatorname{Dom}(\bD_X)\odot
\operatorname{Dom}(\bD_Y)$.

Let $G$ be a finite group acting continuously on $X$.
Following \cite[Def.~2.1]{FR}, we recall
the notion of an equivariant unbounded Kasparov
module defining a class in  $KK^G(A,B)$ in the case of
$A=C_0(X)$, $B=\C$, see
\cite[Remark~2.2]{FR}.

\begin{definition}\label{equivkkdef}\cite{FR}
A $G$-equivariant even
unbounded $K$-cycle defining a class in
\[
K_0^{G,an}(X)=KK_0^G(C_0(X),\C)
\]
consists of a tuple
$(\mathcal H,\Gamma,\rho_X,\bD,U)$ where \((\mathcal H,\Gamma,\rho_X,\mathbf D)\) is an ordinary unbounded Kasparov cycle, together with a unitary representation \(U:G\to\mathcal U(\mathcal H)\) satisfying the following equivariance conditions:

\begin{enumerate}[label=(\alph*)]

  \item There exists a dense $G$-invariant
    $\ast$-subalgebra $A\subset C_0(X)$ such that,
    for every $f\in A$,
    \[
    \rho_X(f)\operatorname{Dom}(\bD)
    \subseteq\operatorname{Dom}(\bD)
    \]
    and $[\bD,\rho_X(f)]$ extends to a bounded
    operator on $\mathcal H$.

    \item  $U$  preserves the grading, i.e.,
    \[
    U_g\Gamma=\Gamma U_g,
    \qquad \forall g\in G
    \]

     \item The representation $\rho_X$ is $G$-equivariant:
    \[
    U_g\rho_X(f)U_g^{-1}
    =\rho_X(g\cdot f)
    \]
    for every $g\in G$ and $f\in C_0(X)$, where
    $(g\cdot f)(x)=f(g^{-1}x)$.

  \item The operator $\bD$ is $G$-equivariant:
    \[
    U_g\operatorname{Dom}(\bD)
    \subset\operatorname{Dom}(\bD),
    \qquad
    U_g\bD U_g^{-1}=\bD,
    \qquad \forall g\in G
    \]
\end{enumerate}
\end{definition}

\subsubsection{Induction of equivariant $K$-cycles}\label{inductkcycles} Let $G$ be a finite group, and $H\subset G$ a subgroup. 

Recall that the induced representation of an $H$-representation $V$  is a $G$-representation defined by $$Ind_H^G(V):=\C[G]\otimes_{\C [H]} V$$ 
we can express
\begin{equation}\label{inductgroup}
    \C[G]\otimes_{\C [H]} V=\bigoplus_{t\in G/H} V_{t}
\end{equation}
where $V_{t}$ is a copy of $V$, and after choosing coset representative $g_tH$ for $t\in G/H$ one can describe the action 
 for $k\in G, v_t\in V_t$ by \cite[Def.~1.12.2]{Sag}
\[k. v_{t}:=((g^{-1}_{kt}kg_t).v)_{kt}  \in V_{kt}\]
note that $g^{-1}_{kt}kg_t\in H$  so the action makes sense.
There is a projection formula
\begin{equation}\label{projectiongroup}
    Ind^G_H(res^G_H(W)\otimes V)=W\otimes Ind^G_H(V)
\end{equation}
for a $G$-representation $W$ and an $H$-representation $V$. Restriction and induction are available on the level of equivariant $KK^G(A,B)$ theory \cite[\S~20.5]{Bl}. Here we only recall the definition in the special case of $A=C_0(X), B=\C$ where $X$ is a $G$-space. We have
\[res^G_H:K_0^{an,G}(X)\rightarrow K_0^{an,H}(X)\]
which can be defined on $K$-cycles by restricting the representation $U_G$ to $U_H$.
 By \cite[Thm~20.5.4]{Bl} there is a homomorphism 
\[\iota_{H\uparrow G}:KK^H(C_0(X),\C)\rightarrow KK^G(Ind^G_{H}(C_0(X)),Ind^G_{H}(\C))\]
we have by \cite[Example~20.5.3 (b)]{Bl} $$Ind^G_H(C_0(X))\cong C_0(G/H)\otimes C_0(X)\cong C_0(G/H\times X),\quad Ind^G_H \C=C(G/H)$$ so we can write
$\iota_{H\uparrow G}:KK^H(C_0(X),\C)\rightarrow KK^G(C_0(G/H\times X),C(G/H))$.
By taking  Kasparov product with the equivariant $K$-homology  fundamental class of the finite \(G\)-set $G/H$
\[[G/H]:=[(l^2(G/H),1,m,0,\lambda)]\in KK^G(C(G/H),\C)\]
where $m$ is the multiplication operator and $\lambda$ is the left permutation representation,
we get a morphism $KK^H(C_0(X),\C)\rightarrow KK^G(C_0(G/H\times X),\C)$. Then the projection to the second factor $G/H\times X\rightarrow X$ induces an equivariant  morphism $$p_*:KK^G(C_0(G/H\times X),\C) \rightarrow KK^G(C_0(X),\C)$$ Therefore we define
\begin{equation}\label{inductform}
    Ind_H^G:=  p_*( \iota_{H\uparrow G} \otimes [G/H]):      K_0^{an,H}(X)\rightarrow K_0^{an,G}(X)
\end{equation}

For the following unbounded-cycle description, assume that
the $H$-equivariant cycle
$(\mathcal H,\Gamma,\rho_X,\bD,U_H)$ admits a dense
$G$-invariant $\ast$-subalgebra $A\subset C_0(X)$ such that
\[
\rho_X(a)\operatorname{Dom}(D)\subseteq\operatorname{Dom}(D),
\qquad
[D,\rho_X(a)]\in\mathcal B(\mathcal H),
\qquad \forall a\in A
\]
where the commutator denotes its bounded extension. We can represent $Ind^G_H[\bD]$ in terms of an equivariant $K$-cycle  by summing the $G/H$-translates of $(\mathcal H,\Gamma,\rho_X,\bD,U_H)$ to obtain a $G$-equivariant cycle as in (\ref{inductgroup}), i.e., in terms of Definition \ref{equivkkdef},

\[Ind^G_H\left((\mathcal{H},\Gamma,\rho_X,\bD,U_H)\right)=( \bigoplus_{t\in G/H} \mathcal{H}_{t},\bigoplus_{t\in G/H}  \Gamma, \rho_{ind},\bigoplus_{t\in G/H}  \bD,Ind^G_H U_H)\]
where $\cH_{t}$ is a copy of $\cH$, and by choosing coset representative $g_tH$ for $t\in G/H$ we have for $v_t\in \cH_t$
\[\rho_{ind}(f)(v_t)=(\rho_X(g_t^{-1}.f)v)_t\in  \cH_t\]
and $\bigoplus_{t\in G/H}  \bD$ acts componentwise on $\bigoplus_{t\in G/H}  Dom(\bD)$. Since by assumption $A$ is $G$-invariant, $g_t^{-1}\cdot a\in A$
for every $a\in A$ and $t\in G/H$. Thus
$\rho_{\mathrm{ind}}(a)$ preserves
$\bigoplus_{t\in G/H}  Dom(\bD)$, and
\[
[\bigoplus_{t\in G/H}  \bD,\rho_{\mathrm{ind}}(a)]
=
\bigoplus_{t\in G/H}
[\bD,\rho(g_t^{-1}\cdot a)]
\]
$Ind^G_H U_H$ is the induced representation
\[Ind^G_H U_H: G\rightarrow \mathcal{U}(Ind^G_H \cH)\cong \mathcal{U}(\bigoplus_{t\in G/H} \mathcal{H}_{t})\]
\[\left(Ind^G_H U_H \right)(k)(v_t)=\left(U_H(g^{-1}_{kt}kg_t)(v)\right)_{kt}\in \cH_{kt}\]

If we begin by a $G$-equivariant Kasparov cycle $(\mathcal{H},\Gamma,\rho_X,\bD,U_G)$ with a decomposition $$\cH=\bigoplus_{i\in I}\cH_i$$ where $I$ is a $G$-set which $G$ acts transitively on $I$ such that $U_G(g)\cH_i=\cH_{gi}$ where $\Gamma, \bD, \rho_X(f)$ for $f\in C_0(X)$ preserve each summand, then we have the orbit-stabilizer description as follows.  Fix $i_0\in I$ and let $H=Stab_G(i_0)$, then there is an isomorphism of $G$-equivariant Kasparov cycles
\begin{equation}\label{indcutfomrulakcycles}
[(\mathcal{H},\Gamma,\rho_X,\bD,U_G)] = Ind^G_H[(\mathcal{H}_{i_0},\Gamma_{i_0},(\rho_X)_{i_0},\bD_{i_0},U_H)]
\end{equation}

We will use the projection formula (\ref{projectiongroup}) on $K$ theory level later in this paper. $K_0^{G(H),an}(X)$ is an $R(G)~\left(R(H)\right)$-module. If $x\in K_0^{H,an}(X), a=[V]\in R(G)$ we have
\begin{equation}\label{projectionfomrk}
    a.Ind^G_H(x)=Ind^G_H(res^G_H(a).x)
\end{equation}

\subsection{Topological, geometric and analytic $K$-homology}
\label{compairsonmaps}
We briefly review  various models for $K$-homology, see \cite{A,ABP1} for a detailed discussion. Let  $X$ be a compact finite CW complex.  There are isomorphisms
\[
 K_0^{top}(X)
\xleftarrow{\underset{\cong}{\phi}}
K_0^{geo}(X)
\xrightarrow{\underset{\cong}{\mu}}
K_0^{an}(X)
\]
A class $[(M,E,f)]\in K_0^{geo}(X)$ in Baum–Douglas K-homology \cite{BD} is represented by a triple where $M$ is a closed even-dimensional $spin^c$ manifold, $E\rightarrow M$ is a complex vector bundle, and $f:M\rightarrow X$ is a continuous map. Equivalence of geometric and analytic $K$-homology \cite{BHS} is given by 
\[\mu:K_0^{geo}(X) \xrightarrow{\cong} K_0^{an}(X) \]
\[[(M,E,f)]\mapsto f_*[\bD_M^E]\] 
where $f_*:K_0^{an}(M) \rightarrow K_0^{an}(X)$ and $[\bD_M^E]$ denotes the  class in $K_0^{an}(M)$ corresponding to the  $spin^c$ operator on $M$ twisted by $E$:
\[[\bD_M^E]:=E\cap[\bD_M^{spin^c}]:=[E]\otimes_{C(M)}[\bD_M^{spin^c}] \in K_0^{an}(M)\]
where, in $KK$-theory, $E$ determines an element
\[
[E]\in KK_0(C(M),C(M))
\]
There is an isomorphism 
\[\phi:K_0^{geo}(X) \xrightarrow{\cong} K_0^{top}(X) \]
\[[(M,E,f)]\mapsto f_*([E]\cap [M]_K)\]
where $\cap$ denotes the $K^0$-module structure of $K^{top}_0$; $\cap:K^* \otimes K_*\rightarrow K_*$, and $[M]_K$ is the fundamental $K$-homology class of $M$.

We denote the Chern character on vector bundles as
\[ch:K^0(X)\rightarrow H^*(X,\mathbb{C})\] and the homological Chern character \cite{BFM2} as
\[ch^{top}:K_0^{top}(X)\rightarrow H_*(X,\mathbb{C})\]
The Chern character on $K_0^{geo}(X)$ is defined by \cite[\S~2.6]{BE}
\[ch^{geo}(M,E,f)=f_*(ch(E)\cup Td^{spin^c}(TM)\cap [M])\]
\[Td^{spin^c}(TM)=e^{c_1(L_{M})/2}\hat{A}(TM)\]
where $L_{M}$ is the determinant line bundle of the $spin^c$ structure. We may define the analytic Chern character by  \cite[\S~2.6]{BE}
\[Ch:=ch^{geo}\circ\mu^{-1}:K_0^{an}(X) \rightarrow  H_*(X,\mathbb{C})\]
 We have the compatibility of topological and analytic Chern characters
\begin{equation}\label{compatiblitcherntodd2}
\begin{tikzcd}
K_0^{geo}(X)
    \arrow[d, "\phi"]
    \arrow[rrr, "\mu"]
    \arrow[drrr, "ch^{geo}"]
&&&
K_0^{an}(X)
    \arrow[d, "Ch"]
\\
K_0^{top}(X)
    \arrow[rrr, "ch^{top}"]
&&&
H_*(X,\mathbb{C})
\end{tikzcd}
\end{equation}
The above comparison maps hold equivariantly when $X$ is a $G$-CW-complex for a finite group $G$, e.g., by \cite{BOSW} 
\begin{equation}\label{equicompare}
    \mu^G:K_0^{G,geo}(X) \xrightarrow{\cong} K_0^{G,an}(X) 
\end{equation}
A class in $K_0^{G,geo}(X)$ is a triple $[(M,E,f)]$ where $M$ is a $G$-$spin^c$ manifold, $E$ is a complex $G$-vector bundle on $M$, and $f:M\rightarrow X$ is $G$-equivariant.

\subsection{Localized characteristic classes}

Let $G$ be a finite group, and let $X$ be a compact $G$-CW-complex. Let $g\in G,\:G'=\langle g\rangle$. Let $X^g$ be the fixed point set of $X$.
\begin{definition}\cite{A} \label{localizedchernclass} The \textit{localized Chern character transformation}

\[Ch(-,g):K_0^{G,an}(X) \rightarrow H_*(X^g,\mathbb{C})\]
is defined as the image under the maps
\[K_0^{G,an}(X)\xrightarrow{res^G_{G'}} K_0^{G',an}(X)\rightarrow K_0^{G',an}(X)_{{g}}\xrightarrow[]{\gamma_{{g}}}K_0^{G',an}(X^g)_{g}\simeq K_0^{an}(X^g) \otimes R(G')_{g}\xrightarrow{Ch\otimes \chi_g}H_*(X^g,\mathbb{C})\]
where $\gamma_g$ is the localization map of analytic $K$-homology \cite{BLP} which is
 the inverse of $$(i^{an})_*:K_0^{G',an}(X^g)_{g} \xrightarrow[]{\cong} K_0^{G',an}(X)_{{g}}$$ induced by the inclusion of the fixed point set.
\[\gamma_{{g}}:=(i^{an})_*^{-1}:K_0^{G',an}(X)_{{g}}\xrightarrow[]{\cong}K_0^{G',an}(X^g)_{g}\]
The \textit{localized Todd transformation}, for a quasi-projective complex variety $X$
\[\td_{*}(-,g):G_0^{G}(X) \rightarrow H^{BM}_*(X^g,\mathbb{C})\]
is defined similarly by localization $L_g$ of equivariant coherent sheaves discussed in  \cite{BFQ,Q79}.
\[G_0^{G}(X)\xrightarrow{res^G_{G'}} G_0^{G'}(X)\rightarrow G_0^{G'}(X)_{{g}}\xrightarrow[]{L_{{g}}}G_0^{G'}(X^g)_{g}\simeq G_0(X^g) \otimes R(G')_{g}\xrightarrow{\td_{*}\otimes \chi_g}  H^{BM}_*(X^g,\mathbb{C})\]
\end{definition}

\subsection{Localization theorem for $spin^c$ manifolds}\label{spinstuff} We start with the case of $spin^c$ manifolds. Later, in the proof of Theorem \ref{problemmoon}, we need a $K$-homology version of the localization theorem for even-dimensional $spin^c$ manifolds.

Let $G=\langle g\rangle$ be a finite cyclic group, $M$ be a closed even-dimensional $G$-manifold equipped with a $G$-equivariant $spin^c$ structure (hence $G$ preserves the $spin^c$ structure). The fixed point set $M^g$ is a submanifold of $M$, assuming it is connected for the sake of notation. 
Let $N=N^M_{M^g}\to M^g$ be the normal bundle.
Assume that $M^g$ is equipped with a $\mathrm{spin}^c$
structure. We give this structure the trivial $G$-lift
and equip $N$ with the complementary $G$-equivariant
$\mathrm{spin}^c$ structure, compatible with the
restriction of the given structure on $M$. Thus (see also \cite[4.4]{BHS})
\begin{equation}\label{spincompatibilty2}
    TM|_{M^g}\cong N\oplus TM^g
\end{equation}
\begin{equation}\label{spincompatibilty}
S_M|_{M^g}\cong S_N\widehat\otimes S_{M^g}
\end{equation}
equivariantly as $\mathbb{Z}_2$-graded Clifford modules,  with $G$ acting trivially on $S_{M^g}$, where $S_M$ ($S_{M^g}$) denotes the spinor bundle of the $spin^c$ bundle $TM$ ($TM^g$), and $S_N=S^{+}_N\oplus S^{-}_N$ is the (equivariant) spinor bundle of the normal bundle $N$, where each $S^{\pm}_N$ is an equivariant complex vector bundle.

By the equivariant $spin^c$ Poincar\'e duality
\cite[Thm. 3.1 and Remark 3.2(1)]{BOSW}, cap product with the
equivariant fundamental Dirac class induces an isomorphism
\begin{equation}\label{equivariant-$spin^c$-Poincar}
    PD_M:K_G^0(M)\xrightarrow{\cong}K_0^{G,\mathrm{an}}(M),
\qquad
[E]\longmapsto [\bD_M^E]
\end{equation}
for a virtual  $G$-equivariant  vector bundle $E$. 
In $KK$-theory, $E$ determines an element
\[
[E]\in KK^G_0(C(M),C(M))
\]
\[
[\bD_M^E]
:=E\cap[\bD_M^{spin^c,G}]:=
[E]\otimes_{C(M)}[\bD_M^{spin^c,G}]
\]
Every element in $K_0^{G,\mathrm{an}}(M)$ is of this form.
Non-equivariantly, the compatibility of Chern characters (\ref{compatiblitcherntodd2})  implies
\begin{equation}\label{chernofPDspin}
    Ch([\bD_M^E])=ch(E)\cup Td^{spin^c}(TM)\cap [M]
\end{equation}
Let $i:M^g \hookrightarrow M$ be the inclusion, and
\[
i^*:K_G^0(M)\rightarrow K_G^0(M^g)
\]
be the equivariant pullback of virtual vector bundles.
\begin{theorem}\label{smoothlocalanal}
With the above assumptions, we have
\[
\gamma_g([\bD_M^E]_g)
=
\left(
\frac{(i^*E\otimes(S^{+}_N-S^{-}_N))_g}
{\lambda_{-1}(N_{\mathbb{C}})_g}
\right)
\cap [\bD^{spin^c}_{M^g}]
\]
\begin{equation}\label{cherncompelxlocalition}
    Ch([\bD_M^E],g)
=
Ch_g\circ\gamma_g([\bD_M^E]_g)
=
\frac{ch_g(i^*E\otimes(S^{+}_N-S^{-}_N))}
{ch_g(\lambda_{-1}(N_{\mathbb{C}}))}
Td^{spin^c}(TM^g)\cap [M^g]
\end{equation}

\end{theorem}
\begin{proof}
We have not found a proof of this formula in the ${spin}^c$ setting in the literature. We therefore derive it from the localization theorem for Riemannian manifolds given in \cite{HW}. Equip $M$ with a $G$-invariant Riemannian metric. By \cite[Prop.~5.20]{HW}, the localized class $[\bD_M^E]_g$ is


\begin{equation}\label{maineq}
    [\bD_M^E]_g=
(j^{TM}_{TM^g})^*_g[\sigma_{\bD_M^E}]_g
\otimes_{C_0(TM^g)}
[\tau^*_{M^g}\Lambda N_{\mathbb{C}}]^{-1}_g
\otimes_{C_0(TM^g)}
[\bD_{TM^g}]_g
\in KK_0^G(C(M),\mathbb{C})_g
\end{equation}
where $[\Lambda N_{\mathbb C}]=[\Lambda^{\mathrm{ev}}N_{\mathbb C}]
-[\Lambda^{\mathrm{odd}}N_{\mathbb C}]=:\lambda_{-1}(N_{\mathbb C})$ which becomes invertible after localization at $g$, 
\[
\tau_{M^g}: TM^g\rightarrow M^g
\]
\[
[\bD_{TM^g}] \in KK_0^G(C_0(TM^g),\mathbb{C})
\]
$\bD_{TM}$ of a Riemannian manifold $M$ is the Dolbeault--Dirac operator on $TM$, and
\[
(j^{TM}_{TM^g})^*:
KK_0^G(C(M),C_0(TM))
\rightarrow
KK_0^G(C(M),C_0(TM^g))
\]
The principal symbol $\sigma_{\bD_M^E}$ of $\bD_M^E$ defines a class 
\[[\sigma_{\bD_M^E}] \in KK_0^G(C(M),C_0(TM))\]
 We have
 \begin{equation}\label{eqinproofnedd}
     (j^{TM}_{TM^g})^*_g[\sigma_{\bD_M^E}]_g
=
[i^*]\otimes[\sigma_{\bD_M^E}|_{TM^g}]
 \end{equation}
where
\[
[i^*] \in KK_0^G(C(M),C(M^g))
\]
is the class associated with the homomorphism of $C^*$-algebras
$C(M)\rightarrow C(M^g)$ induced by the restriction, and
\[
[\sigma_{\bD_M^E}|_{TM^g}]
\in
KK_0^G(C(M^g),C_0(TM^g))
\]
By applying $\gamma_g=(i_*)^{-1}$ to (\ref{maineq}), where $
i_*(-)=[i^*]\otimes -
$, and using (\ref{eqinproofnedd})
we get
\begin{equation}\label{eq2}
    \gamma_g[\bD_M^E]_g=
[\sigma_{\bD_M^E}|_{TM^g}]_g
\otimes_{C_0(TM^g)}
[\tau^*_{M^g}\Lambda N_{\mathbb{C}}]^{-1}_g
\otimes_{C_0(TM^g)}
[\bD_{TM^g}]_g
\in KK_0^G(C(M^g),\mathbb{C})_g
\end{equation}
Since $[\bD_M^E]=E\cap [\bD^{spin^c,G}_M]$,
its symbol is
$
\sigma_{\bD_M^E}
=
[E]\otimes \sigma_{\bD^{spin^c,G}_M}
$, therefore
\[
[\sigma_{\bD_M^E}|_{TM^g}]_g
=
[(i^*E)_g]
\otimes_{C(M^g)}
[\sigma_{\bD^{spin^c,G}_M}|_{TM^g}]_g
\]
By the compatibility of equivariant $spin^c$ structures (\ref{spincompatibilty2}) and (\ref{spincompatibilty}), we have (by multiplicativity of the $spin^c$ symbol under direct sums)
\[
[\sigma_{\bD^{spin^c,G}_M}|_{TM^g}]_g
=
[S^+_N-S^-_N]_g
\otimes_{C(M^g)}
[\sigma_{\bD^{spin^c}_{M^g}}]_g
\in
KK_0^G(C(M^g),C_0(TM^g))
\]
and therefore
\[
[\sigma_{\bD_M^E}|_{TM^g}]_g
=
[i^*E]_g
\otimes_{C(M^g)}
[S^+_N-S^-_N]_g
\otimes_{C(M^g)}
[\sigma_{\bD^{spin^c}_{M^g}}]_g
\in
KK_0^G(C(M^g),C_0(TM^g))
\]
Associativity of the Kasparov product and the following projection formula (which means twisting the symbol by a vector bundle on $X^g$ is equivalent to twisting it by the pullback of that bundle to $TX^g$)
\[
[\sigma_{\bD^{spin^c}_{M^g}}]\otimes_{C_0(TM^g)}[\tau^*_{M^g}\Lambda N_{\mathbb{C}}]^{-1}_g
=
\lambda_{-1}(N_{\mathbb{C}})_g^{-1}\otimes_{C(M^g)}[\sigma_{\bD^{spin^c}_{M^g}}]
\]
imply
\[
\gamma_g[\bD_M^E]_g
=
\left(
[i^*E]_g
\otimes_{C(M^g)}
[S^+_N-S^-_N]_g
\otimes_{C(M^g)}
\lambda_{-1}(N_{\mathbb{C}})_g^{-1}
\right)
\otimes_{C(M^g)}
[\sigma_{\bD^{spin^c}_{M^g}}]_g
\otimes_{C_0(TM^g)}
[\bD_{TM^g}]_g
\]
We apply a Kasparov index theorem \cite[\S4]{HW} on $M^g$, which says
\[
[\sigma_{\bD^{spin^c}_{M^g}}]
\otimes_{C_0(TM^g)}
[\bD_{TM^g}]
=
[\bD^{spin^c}_{M^g}]
\in
KK_0^G(C(M^g),\mathbb{C})
\]
Combining everything, we get
\[
\gamma_g[\bD_M^E]_g
=
\left(
[i^*E]_g
\otimes_{C(M^g)}
[S^+_N-S^-_N]_g
\otimes_{C(M^g)}
\lambda_{-1}(N_{\mathbb{C}})_g^{-1}
\right)
\cap
[\bD^{spin^c}_{M^g}]
\]
(\ref{cherncompelxlocalition}) then follows from the above formula and (\ref{chernofPDspin}).
\end{proof}

\subsection{Thom-Mather stratified spaces and smoothly stratified spaces}

We briefly recall the notions of Thom--Mather stratified spaces
and smoothly stratified spaces, following \cite{ABP1}, and $G$-Thom--Mather stratified space in the sense of \cite{BLP}.

A \emph{Thom--Mather stratified space} is a triple
$(X,\mathcal S,\mathcal T)$, where $X$ is a locally compact,
second-countable Hausdorff space, $\mathcal S$ is a locally
finite partition of $X$ into locally closed smooth manifolds,
called strata, satisfying the frontier condition
\[
S\cap\overline{T}\neq\emptyset
\quad\Longrightarrow\quad S\subseteq\overline{T}
\]
and $\mathcal T$ is a collection of compatible control data.
More precisely, for each stratum $S$, the control data consist
of a tubular neighborhood $T_S$ of $S$, a continuous retraction
$\pi_S:T_S\to S$, and a continuous function
$\rho_S:T_S\to[0,\infty)$ satisfying
\[
\pi_S|_S=\mathrm{id}_S,\qquad \rho_S^{-1}(0)=S
\]
If $S\subseteq\overline{T}$ and $S\neq T$, the restriction
\[
(\pi_S,\rho_S):T_S\cap T\longrightarrow
S\times(0,\infty)
\]
is required to be a smooth submersion. Furthermore, for
$S\subseteq\overline{T}$ and $T\subseteq\overline{U}$,
the control data satisfy the compatibility relations
\[
\pi_S\circ\pi_T=\pi_S,\qquad
\rho_S\circ\pi_T=\rho_S
\]
whenever the compositions are defined on the relevant
intersections.

A \emph{ smoothly stratified space} or \emph{smooth atlas stratified space}, as introduced in
\cite[Definition 2.2]{ABP1}, is a conically stratified space
equipped with an atlas of compatible stratified charts of
the form
\[
\mathcal U\cong\mathbb R^k\times C(L)
\]
where $C(L)$ denotes the open cone on a compact stratified
link $L$ of smaller depth. The smooth structure is defined
inductively on the depth, requiring the transition maps
and their inverses to admit smooth stratified lifts after
unfolding the conical variables.

By \cite[Theorem 2.13]{ABP1}, Thom--Mather stratified spaces
and smoothly
stratified spaces  describe the same class
of spaces.

An important advantage of the smooth
atlas description is that this class is closed under
Cartesian products. Indeed, this follows from the local
identification
\[
C(L_1)\times C(L_2)\cong C(L_1*L_2)
\]
where $L_1*L_2$ denotes the stratified join. In particular, if $X$ is smoothly stratified, then so is $X^n$ for every
$n\geq 1$. In contrast, while the  product of Thom-Mather stratified spaces is
again a Thom-Mather stratified space, the product of the Thom-Mather data is \emph{not} Thom-Mather data for the
product as noted in \cite{Ve} and \cite[Remark 2.12]{ABP1}.

A smoothly stratified space $X$ is a \textit{pseudomanifold} if its top-dimensional stratum is dense in $X$ and every other stratum has codimension at least $2$, and  it is furthermore called a
\emph{Witt space} if, for every stratum $S$ of odd
codimension $2k+1$, its link $L_S$ satisfies
\[
IH_k^{\bar m}(L_S;\mathbb Q)=0
\]
where $\bar m$ denotes the lower-middle perversity
and $IH_*^{\bar m}$ denotes intersection homology
with rational coefficients. In particular if all strata of $X$ are even dimensional, such as complex algebraic varieties, then $X$ is a Witt space. The Witt condition does not depend on the stratification. A product of Witt spaces is again Witt.

Now assume we have a finite group $G$. Following \cite[Definition 2.1]{BLP}, a
\emph{Thom--Mather $G$-space} is a Thom--Mather stratified space
$(X,\mathcal S,\mathcal T)$ equipped with a continuous
$G$-action such that every $g\in G$ preserves each
stratum $S\in\mathcal S$ and restricts to a diffeomorphism
of $S$. Moreover, the control data are $G$-equivariant:
for every $S\in\mathcal S$ and $g\in G$,
\[
g(T_S)=T_S,\qquad
\pi_S(gx)=g\pi_S(x),\qquad
\rho_S(gx)=\rho_S(x)
\]

The symmetric group $S_n$ acts on $X^n$ by permutation of coordinates. We need the following fact in this paper.

\begin{theorem}\label{Sn-mather-X^n}
Let $X$ be an even-dimensional, oriented, closed, smoothly stratified Witt space. Then $X^n$ is a Thom--Mather Witt space and admits an $S_n$-Thom--Mather structure (in the sense of \cite{BLP}). Furthermore, $X^n$ admits an $S_n$-equivariant triangulation.
\end{theorem}

\begin{proof}
The product of closed smoothly stratified Witt spaces is again such a space; see, for example, \cite[Corollary~5.18]{ABP1}. The strata of $X^n$ are of the form
$S_{\alpha}:=S_{\alpha_1}\times\ldots\times S_{\alpha_n}$,
where each $S_{\alpha_i}$ is a stratum of $X$. The group $S_n$ acts by permuting the factors.
We need to show that $X^n$ admits an $S_n$-Thom--Mather structure. 
Thom--Mather spaces are equivalent to Whitney-embedded stratified spaces \cite[\S 2.4]{ABP1}. Thus, we may assume that $X$ is a Whitney-stratified subset of a smooth manifold $M$, with an embedding $X\hookrightarrow M$. Then $X^n\hookrightarrow M^n$ is an $S_n$-equivariant embedding.
An element $\sigma\in S_n$ may not preserve each stratum. For this reason, we consider the $S_n$-orbits of the product strata:
\[
\tilde{S}_{\alpha}:=\bigcup_{\sigma \in S_n} S_{\alpha_{\sigma^{-1}(1)}}\times \ldots\times S_{\alpha_{\sigma^{-1}(n)}}
\]
Then $\{\tilde{S}_{\alpha}\}$ is a Whitney stratification in which each stratum is preserved by every element of $S_n$. Then \cite[Thm.~2.2]{BLP} and \cite[Thm.~2.12]{PW} imply that $X^n$ admits $S_n$-equivariant control data, making it a Thom--Mather $S_n$-stratified space.
It remains to establish the existence of an
$S_n$-equivariant triangulation. By Goresky's
triangulation theorem \cite{G} (see also
\cite[Proposition~4.5]{ABP1}), $X$ admits a
finite triangulation $K$ compatible with its
stratification. The product $|K|^n$ is a compact semialgebraic
set, and the natural $S_n$-action by permutation
of coordinates is semialgebraic and proper.
Hence, by \cite[Prop.~6.7]{BLP},
$|K|^n$ admits an $S_n$-equivariant
triangulation.
Finally, since $\dim X$ is even, permutations
of the factors preserve the product orientation.
\end{proof}

\subsection{Equivariant $K$-homology class of signature operator and equivariant  $L$-classes}

Let $X$ be a compact oriented even-dimensional  smoothly stratified Witt space. In particular, $K_{\dim X}^{an}(X)\simeq K_0^{an}(X)$.
\begin{theorem}
    \cite{ALMP} There is a $K$-homology class of signature operator $$\left[D_X^{sign}\right]\in K^{an}_0(X)$$ by choosing a wedge metric. The class does not depend on the metric.
\end{theorem}
More precisely, in terms of Definition
\ref{kkdef}, we consider the tuple
\begin{equation}\label{singaureclass}
    (\mathcal H,\Gamma,\rho_X,\bD)
\end{equation}
where
$
\mathcal H
=L^2\Omega^*(X_{\mathrm{reg}}, g),
$
for a   wedge metric $g$ on $X_{reg}$. The grading $\Gamma$ is the signature
grading $\tau$ introduced in (\ref{siggradingdef}). The representation
$
\rho_X:C(X)\longrightarrow\mathcal B(\mathcal H)
$
is given by multiplication:
$
\rho_X(f)\omega=f\omega.
$
 $\bD=D^{sign}_{\mathbf g}$ is
the signature operator associated with
the de Rham operator $d+d^*$, graded by
$\tau$, with the self-adjoint VAPS domain
\[
\operatorname{Dom}(\bD)
=\mathcal D_{VAPS}(D^{sign}_{\mathbf g})
\]
The Witt condition implies that $\bD$
is self-adjoint and has compact resolvent.

Now let $G$ be a finite group and assume  $X$ is a compact oriented even-dimensional Thom-Mather G-Witt space where $G$ preserves the orientation, in the sense of \cite[\S2]{BLP}. 
\begin{theorem}\label{equivsingaureclass}
    \cite[Thm. 4.1]{BLP} There is an equivariant signature $K$-homology class $$\left[D_X^{sign,G}\right]\in K^{G,an}_0(X)$$ by choosing a $G$-invariant wedge metric. The class does not depend on the metric.
\end{theorem}
The extra data other than those mentioned in (\ref{singaureclass}) needed to define an equivariant class in terms of Definition
\ref{equivkkdef} is as follows.
First, we choose a  $G$-invariant wedge metric $g=\mathbf g$ on $X_{reg}$. 
The action of $G$ on $X_{\mathrm{reg}}$
induces a unitary representation
\begin{equation}\label{equivrepdef}
    U:G\longrightarrow\mathcal U(\mathcal H),
\qquad
U_g\omega=(g^{-1})^*\omega
\end{equation}
Since $\mathbf g$ is $G$-invariant and
the action preserves the orientation,
$U_g$ commutes with the signature grading
$\Gamma$. Moreover,
$
U_g\rho_X(f)U_g^{-1}
=\rho_X(g\cdot f)
$
where $(g\cdot f)(x)=f(g^{-1}x)$. Furthermore, since $\mathbf g$ is
$G$-invariant, the signature operator
and its VAPS domain are $G$-invariant:
\[
U_g\operatorname{Dom}(\bD)
=\operatorname{Dom}(\bD),
\qquad
U_g\bD U_g^{-1}=\bD
\]
 the associated class is represented by
\[
[D_X^{sign,G}]
\in KK_0^G(C(X),\C)=K^{G,an}_0(X)
\]

\begin{theorem}\label{Banagl_localization_comparison}\cite{A}
Let $G$ be a finite group and let $X$ be a compact oriented even-dimensional triangulated Witt pseudomanifold on which $G$ acts by orientation-preserving simplicial maps. Let $g\in G$, and 
$
j:X^g\rightarrow X
$
be the inclusion of the fixed point set. Then
\[
L^{}_*(g,X)=j_*(\Psi_2(L^{an}_*(X,g)))
\]
where $L_*(g,X)$ is the Banagl-Zagier equivariant $L$ class \cite{B} and
\[
L^{an}_*(X,g):=Ch(\left[D_X^{sign,G}\right],g)
\in H_*(X^g,\mathbb{C})
\]
\end{theorem}
Banagl \cite{B} showed that $X/G$ is also a Witt pseudomanifold. By applying Theorem \ref{Banagl_localization_comparison} to \cite[Thm. 10.1]{B}  we get:
  \begin{corollary}\label{averageLclass}\cite[Corollary 4.6]{A}
  Let $G$ be a finite group and let $X$ be a compact oriented even-dimensional triangulated Witt pseudomanifold on which $G$ acts by orientation-preserving simplicial maps. Then
    \begin{equation}\label{averageLclasseq}
        L^{an}_*(X/G)=\frac{1}{|G|}\sum_{g\in G} \pi^g_*Ch(\left[D_X^{sign,G}\right],g)=\frac{1}{|G|}\sum_{g\in G} \pi^g_*(L^{an}_*(X,g))
    \end{equation}
    where $\pi^g:X^g\hookrightarrow X\rightarrow X/G$ is the projection.
\end{corollary}

\section{Power operations in $K$-homology}\label{poweropsection}

We recall power operations in $K$-theory and $K$-homology. We start with the $K$-theory group $K^0(X)$.  Atiyah's power operation \cite{At}
\[P_n:K^0(X)\rightarrow K_{S_n}^0(X^n)\]
can be described by the external product in $K$-theory
\[\boxtimes : K^0(X)\otimes K^{0}(Y)\rightarrow K^0{(X\times Y)} \]
We identify $$K^0_{S_n}(X^n)\cong K_0(D_b(Vec^{S_n}(X^n)))$$ i.e., the Grothendieck group of the bounded derived category of complexes of $S_n$-equivariant vector bundles, e.g., see \cite[\S6]{CMSS}. If $[E]\in K^0(X)$ is an honest vector bundle, then $P_n(E):=[E^{\boxtimes n}]$ and
a permutation $\sigma\in S_n$ acts on it by permuting the factors without sign. If    $E^{\bullet}$ is a bounded  complex of vector bundles with the corresponding class $$[E^{\bullet}]=\sum_i (-1)^i[E^i]\in K^{0}(X)$$ then
\[P_n([E^{\bullet}]):=[(E^{\bullet})^{\boxtimes n}]\]
where $(E^{\bullet})^{\boxtimes n}$ is the total external tensor product complex with
\[((E^{\bullet})^{\boxtimes n})^k=\oplus_{i_1+\ldots+i_n=k} E^{i_1}\boxtimes \ldots \boxtimes E^{i_n}\]
and the differential is the usual total complex differentials (with sign). $S_n$ acts on $(E^{\bullet})^{\boxtimes n}$ by the graded permutation; for $e_j\in E^{i_j}, \sigma\in S_n$
\[\sigma(e_1\boxtimes\ldots \boxtimes e_n)=\epsilon(\sigma,i_1,\ldots,i_n)e_{\sigma^{-1}(1)}\boxtimes \ldots \boxtimes e_{\sigma^{-1}(n)}\]
where $\epsilon(\sigma,i_1,\ldots,i_n)=(-1)^{\sum_{a<b,\sigma(a)>\sigma(b)}i_ai_b}$.  $(E^{\bullet})^{\boxtimes n}$ is an $S_n$-equivariant complex of vector bundles on $X^n$ defining an element $[(E^{\bullet})^{\boxtimes n}]\in K_{S_n}^0(X^n)$. Note that the power operation is not additive.\\

By a similar process of defining external products of bounded coherent complexes, we can define power operations on $G_0(X)=K_0(Coh(X))$  
\begin{equation}\label{powerop_G0}
    P_n:G_0(X)\rightarrow G^{S_n}_0(X^n)
\end{equation}
by identifying $G^{S_n}_0(X^n)\cong K_0(Coh^{S_n}(X^n))$ and from the external product on $G_0$ \cite[\S3.1]{BFM2}
\[\boxtimes: G_0(X)\otimes G_0(Y)\rightarrow G_0(X\times Y)\]

Now we describe the power operation in analytic $K$-homology.   Consider 
the $n$-fold external product (see (\ref{external_kasparov_product})); given $[x]=[\cH,\Gamma,\rho_X,\bD] \in K^{an}_0(X)$, 
\begin{equation}\label{poweropanal}
    [x]^{\boxtimes n}=[\cH^{\widehat{\otimes}n},\Gamma^{\otimes n},\rho_X^{\widehat{\otimes}n},\bD_{\Gamma}^{(n)},U_{\Gamma}^n] \in K^{S_n,an}_0(X^n)
\end{equation}
where
  $$D_{\Gamma}^{(n)}:=\sum^n_{j=1} 1\widehat{\otimes} \ldots \widehat{\otimes} D\widehat{\otimes} \ldots \widehat{\otimes} 1=\sum_{j=1}^n
\Gamma^{\otimes(j-1)}
\otimes \bD \otimes 1^{\otimes(n-j)}$$
and  the $S_n$-equivariant structure
\[
U_{\Gamma}^n:S_n\to
    \mathcal U(\mathcal{H}_X^{\widehat{\otimes}n})
    \]
    is given, for elements
    $\xi_i\in\mathcal H$ of grade
    $|\xi_i|_{\Gamma}\in\mathbb Z_2$, by
    \begin{equation}\label{gpsign}
         U_{\Gamma,\sigma}^{n}(\xi_1\widehat{\otimes}\cdots
    \widehat{\otimes}\xi_n)
    =
    (-1)^{\Gamma(\sigma;\xi)}
    \xi_{\sigma^{-1}(1)}\widehat{\otimes}
    \cdots\widehat{\otimes}
    \xi_{\sigma^{-1}(n)}
    \end{equation}
    where
    \[
    \Gamma(\sigma;\xi)
    =
    \sum_{\substack{i<j\\
    \sigma(i)>\sigma(j)}}
    |\xi_i|_{\Gamma}|\xi_j|_{\Gamma}
    \quad(\mathrm{mod}\ 2)
    \]
    The required conditions in Definition \ref{equivkkdef} such as $U_{\Gamma}^n$ commuting with $\Gamma^{\otimes n}$, $\rho^{\widehat{\otimes}n}$ being $S_n$-equivariant, and $\bD_{\Gamma}^{(n)}$ being $S_n$-equivariant are clear, hence it defines an element in $K^{S_n,an}_0(X^n)$.

 Therefore, we get (see \cite[Lemma~2.4]{CEKN} for a general statement)
\begin{equation}\label{powerop_ank}
    P_n:K^{an}_0(X)\rightarrow K_0^{S_n,an}(X^n)
\end{equation}
\[[x]\rightarrow [x]^{\boxtimes n}\]

We also describe power operation in the language of Baum–Douglas K-homology which shows the connection to the power operation in $K$-theory. The external product
 \begin{equation}\label{external_kasparov_product_baumdouglass}
     \boxtimes:K^{geo}_0(X)\otimes K^{geo}_0(Y)\rightarrow K^{geo}_{0}(X\times Y)
 \end{equation}
 is defined by (e.g., see \cite[\S~3.1]{J} for a general definition)
 \[[(M_1,E_1,f_1)]\boxtimes [(M_2,E_2,f_2)]=[(M_1\times M_2,E_1\boxtimes E_2,f_1\times f_2)]\]
 where $E_1\boxtimes E_2$ is the external product of vector bundles, and $M_1\times M_2$ carries the product $spin^c$-structure.  Given  $[(M,E,f)]\in K_0^{geo}(X)$
 \[[(M,E,f)]^{\boxtimes n}=[(M^n,E^{\boxtimes n},f^n)]\in K_0^{S_n,geo}(X^n)\]
 Indeed, $S_n$ permutes the factors of $X^n$, $f^n$ is $S_n$-equivariant, and since $M$ is even dimensional the product $spin^c$-structure is preserved. Therefore this defines a class in $K_0^{S_n,geo}(X^n)$. Under the isomorphism
\begin{equation}\label{an_bdmapinsec7}
    \mu:K_0^{geo}(X) \xrightarrow{\cong} K_0^{an}(X)
\end{equation} 
 (Section \ref{compairsonmaps}), for a given   $\mathcal{F}\in K_0^{an}(X)$,  there is a class 
 \begin{equation}\label{f=g,bm=kk}
   [(M,E,f)]\in K_0^{geo}(X)\qquad  \mathcal{F}=\mu([(M,E,f)])=f_*(\mathcal{G})
 \end{equation}
where $f_*:K_0^{an}(M) \rightarrow K_0^{an}(X)$ and
$$\mathcal{G}:=[\bD_M^E]=E\cap[\bD_M^{spin^c}]=[E]\otimes_{C(M)}[\bD_M^{spin^c}] \in K_0^{an}(M)$$ Note that $[\bD_M^E]\in K_0^{an}(M)$ is the class obtained from  $[(M,E,id)]$  under the isomorphism $$ K_0^{geo}(M) \xrightarrow{\cong} K_0^{an}(M)$$
 As mentioned before this  holds for any compact CW-complex $X$. 
 
 The isomorphism (\ref{an_bdmapinsec7}) is compatible with external products. Indeed, the Kasparov product of fundamental classes associated to $spin^c$-structures is the fundamental class associated to the product $spin^c$-structure \cite[Prop.~11.2.13]{HR},
 \[[\bD^{spin^c}_{M_1}]\boxtimes [\bD^{spin^c}_{M_2}]=[\bD^{spin^c}_{M_1\times M_2}]\]
 again we recall we are using  $K_{even}\cong K_0$ to write all the classes in $K_0$, since $M_1,M_2$ are even-dimensional. Therefore
\begin{align*}
    [\bD^{E_1}_{M_1}]\boxtimes [\bD^{E_2}_{M_2}]=&([E_1]\otimes_{C(M_1)}[\bD_{M_1}^{spin^c}]) \boxtimes ([E_2]\otimes_{C(M_2)}[\bD_{M_2}^{spin^c}])\\
    =&([E_1]\boxtimes[E_2] )\otimes_{C(M_1\times M_2)} ([\bD^{spin^c}_{M_1}]\boxtimes [\bD^{spin^c}_{M_2}])\\
    =& [E_1\boxtimes E_2] \otimes_{C(M_1\times M_2)} [\bD^{spin^c}_{M_1\times M_2}]\\
    =&[\bD^{E_1\boxtimes E_2}_{M_1\times M_2}]
\end{align*}
 Under the equivariant isomorphism (\ref{equicompare})
 \[\mu^{S_n}:K_0^{S_n,geo}(X^n) \xrightarrow{\cong} K_0^{S_n,an}(X^n) \]
\[\mu^{S_n}([(M,E,f)]^{\boxtimes n})=(f^n)_*(\mathcal{G}^{\boxtimes n})=(f^n)_*(E^{\boxtimes n}\cap[D_{M^n}^{spin^c,S_n}])\]

\subsection{Proof of Theorem \ref{analyticmoonen}}

We first prove the theorem for a  closed even-dimensional $spin^c$ manifold $M$ and $\mathcal{G}=[\bD_M^E]=E\cap[\bD_M^{spin^c}]$. The manifold $M^n$ is naturally an orientation-preserving
$S_n$-equivariant $spin^c$ manifold; the $spin^c$ structure on $M^n$
is the product structure, $S_n$ acts on the product spinor bundle by
graded permutation, and the even-dimensionality of $M$ implies that
the permutation preserves orientation. Then the diagonal map
\[
i:\Delta(M)=(M^n)^{\sigma_n}\hookrightarrow M^n
\]
is a smooth equivariant embedding. We can  view the equivariant class $ Ch(-,\sigma_n)$  as $$K_0^{S_n,an}(M^n) \xrightarrow{res_{G'}^{S_n}}  K_0^{G',an}(M^n) \rightarrow K_0^{G',an}(M^n)_{{g}}\xrightarrow{\gamma_{g}}  K_0^{G',an}(\Delta(M))_{{g}} \cong  K_0^{an}(M)\otimes R(G')_{{g}}\xrightarrow{Ch_g} H_*(M,\mathbb{C})$$
where $g=\sigma_n,\: G'=\langle \sigma_n\rangle$.
\begin{theorem}\label{problemmoon} Let  $M$ be a closed even-dimensional $spin^c$ manifold, and  let $ \sigma_n\in S_n$ be an $n$-cycle. Identify $M\cong \Delta(M)\subset M^n$. Let $E$ be a complex vector bundle on $M$ and consider 
$\mathcal{G}=E\cap[\bD_M^{spin^c}]:=[E]\otimes_{C(M)}[\bD_M^{spin^c}]\in K_0^{an} (M)$. Then,   the following equality holds in $H_*(M,\mathbb{C})$ 
\begin{equation}
 Ch(P_n(\mathcal{G}),\sigma_n)=\Psi_n Ch(\mathcal{G})
\end{equation} 
\end{theorem}
\begin{proof}
Let $N=N^{M^n}_{\Delta(M)}$ be the normal bundle of $i:(M^n)^{\sigma_n}=\Delta(M)\rightarrow M^n$. The
$spin^c$ structure on $T(M^n)|_{\Delta(M)}$ and $T\Delta(M)\cong TM$ induce a
complementary $spin^c$ structure on $N$.  We equip the diagonal $TM$-factor with the trivial $G'$-lift; the complementary
$spin^c$ structure on $N$ therefore contains all the nontrivial
graded-permutation action. 
\begin{equation}\label{copatibliltynormalproduct}
N\oplus TM
\cong
i^*T(M^n)
\cong
TM^{\oplus n}
\cong
TM\otimes\rho_n
\end{equation}
\begin{equation}\label{N_of_symmetric}
N\cong TM\otimes(\rho_n-1)
\end{equation}
where
\[
\rho_n:G'\longrightarrow O(\mathbb{R}^n)
\]
is the permutation representation of
$G'=\langle\sigma_n\rangle$, and $1$ denotes the trivial
representation.  Recall
\[\mathcal{G}^{\boxtimes n}=E^{\boxtimes n}\cap[D_{M^n}^{spin^c}]\]
By Theorem
\ref{smoothlocalanal}, applied to $M^n$ with
$(M^n)^g=\Delta(M)\cong M$, we get
\[
Ch(P_n(\mathcal{G}),\sigma_n)
=
\left(
ch_g(i^*E^{\boxtimes n})
\frac{ch_g(S_N^+-S_N^-)}
{ch_g(\lambda_{-1}(N_{\mathbb{C}}))}
Td^{spin^c}(TM)
\right)\cap[M]
\]
where
$
S_N=S_N^+\oplus S_N^-
$
is the equivariant spinor bundle of the normal bundle. By \cite[Eq.~(3.1)]{At}, $i^*E^{\boxtimes n} \cong E^{\otimes n}$. By  \cite[p.~179]{At}, under the evaluation map ($g=\sigma_n$)
\[\chi_g:K_{S_n}^0(\Delta(M))\cong K^0(M)\otimes R(S_n)\rightarrow K^0(M)\otimes \mathbb{C} \]
\[\chi_g(E^{\otimes n})= \psi^n E\]
where $\psi^n$ is the $n$-th Adams operation in $K^0(M)$ \cite[Prop.~3.2]{At}. By \cite[\S~7]{At}, $ch(\psi^n E)
=
\Psi^n(ch(E))$. Putting together everything we get 
\[ch_g(i^*E^{\boxtimes n})=\Psi^n(ch(E))\]
Say
$
\dim_{\mathbb{R}}M=2d
$
with formal Chern roots of $TM_{\mathbb{C}}$ given by
$\{\pm p_a\}_{a=1}^d$. By Lemma \ref{chofeulerfactor} we have
\begin{equation}
\frac{ch_g(S_N^+-S_N^-)}
{ch_g(\lambda_{-1}(N_{\mathbb{C}}))}
=
e^{(n-1)c_1(L_M)/2}
\prod_{a=1}^d
\frac{e^{p_a/2}-e^{-p_a/2}}
{e^{np_a/2}-e^{-np_a/2}}
\end{equation}
where $L_M$ is the determinant line bundle of the $spin^c$ structure
on $M$. Also
\[
Td^{spin^c}(TM)
=
e^{c_1(L_M)/2}\hat{A}(TM)
=
e^{c_1(L_M)/2}
\prod_{a=1}^d
\frac{p_a}{e^{p_a/2}-e^{-p_a/2}}
\]
Hence
\[
\frac{ch_g(S_N^+-S_N^-)}
{ch_g(\lambda_{-1}(N_{\mathbb{C}}))}
Td^{spin^c}(TM)
=
e^{nc_1(L_M)/2}
\prod_{a=1}^d
\frac{p_a}{e^{np_a/2}-e^{-np_a/2}}
\]
\begin{equation}\label{chgcomputation}
Ch(P_n(\mathcal{G}),\sigma_n)
=
\left(
\Psi^n(ch(E))
e^{nc_1(L_M)/2}
\prod_{a=1}^d
\frac{p_a}{e^{np_a/2}-e^{-np_a/2}}
\right)\cap[M]
\end{equation}
On the other hand, for $a\in H^*(M,\mathbb{C})$, we have (see Section \ref{backgroundnotations} notations)
\[
\Psi_n(a\cap[M])
=
n^{-d}\Psi^n(a)\cap[M]
\]
where $\Psi^n$ is the cohomological Adams operation. Also
\[
\Psi^n(Td^{spin^c}(TM))
=
\Psi^n\left(
e^{c_1(L_M)/2}\hat{A}(TM)
\right)
=
e^{nc_1(L_M)/2}
\prod_{a=1}^d
\frac{np_a}{e^{np_a/2}-e^{-np_a/2}}
\]
By 
(\ref{chernofPDspin}), 
\[
Ch(\mathcal G)
=
ch(E)Td^{spin^c}(TM)\cap[M]
\]
Therefore
\[
\begin{aligned}
\Psi_n Ch(\mathcal G)
&=
n^{-d}\Psi^n
\left(
ch(E)Td^{spin^c}(TM)
\right)\cap[M]
\\
&=
n^{-d}
\left(
\Psi^n(ch(E))
e^{nc_1(L_M)/2}
\prod_{a=1}^d
\frac{np_a}{e^{np_a/2}-e^{-np_a/2}}
\right)\cap[M]
\\
&=
\left(
\Psi^n(ch(E))
e^{nc_1(L_M)/2}
\prod_{a=1}^d
\frac{p_a}{e^{np_a/2}-e^{-np_a/2}}
\right)\cap[M]\\
&=Ch(P_n(\mathcal{G}),\sigma_n)
\end{aligned}
\]
\end{proof}
\begin{lemma}\label{chofeulerfactor} With the notation used in the proof of Theorem \ref{problemmoon}, we have 
    \[\frac{ch_g(S^{+}_N-S^{-}_N)}{ch_g(\lambda_{-1}(N_{\mathbb{C}}))}=e^{(n-1)c_1(L_M)/2}\prod^d_{a=1} \frac{e^{p_a/2}-e^{-p_a/2}}{ (e^{np_a/2}-e^{-np_a/2})}\]
\end{lemma}
\begin{proof}
Let $i:M\cong \Delta(M) \rightarrow M^n$. As equivariant $spin^c$ bundles on $\Delta(M)$, by (\ref{copatibliltynormalproduct}), we have
\[
 N\oplus TM \cong i^*T(M^{ n})\cong TM^{\oplus n}  \cong TM\otimes \rho_n
\]
$N$ carries the \(G'=\langle\sigma_n\rangle\)-equivariant $spin^c$ structure induced by the product $spin^c$ structure on $M^n$, and has determinant line bundle $L_N$ satisfying
\[
L_N\otimes L_M \cong L_M^{\otimes n}
\]
Hence as ordinary complex line bundles,
\[
L_N\cong L_M^{\otimes (n-1)},
\qquad
c_1(L_N)=(n-1)c_1(L_M)
\]
Let $S$ be a notation for spinor bundles. Let $\epsilon_M:=[S^+_{M}]-[S_{M}^-],
\: \epsilon_N:=[S^+_N]-[S_N^-] \in K_{G'}^0(M)$, where $\epsilon_M$ has the trivial $G'$ action. Since $M^n$ has the product $spin^c$ structure, $\epsilon_{M^n}=\epsilon_M^{\boxtimes n}$.
Multiplicativity of spinor modules for orthogonal direct sums gives
\[ \epsilon_N \otimes  \epsilon_M =i^*(\epsilon_M^{\boxtimes n})\cong \epsilon_M^{\otimes n}\]
where $g=\sigma_n$ acts by the graded cyclic permutation of the spinor factors. By \cite{At}
\[ch_g(\epsilon_M)ch_g(\epsilon_N)=ch_g(\epsilon_M^{\otimes n})=ch(\psi^n(\epsilon_M))\]
For a real rank-$2d$ $spin^c$ vector bundle $TM$ with spinor bundle $S_M=S_M^+\oplus S^-_M$, determinant line bundle $L_M$, and formal Chern roots $\{\pm p_a\}^d_{a=1}$ of $TM_{\mathbb{C}}$,
we have \cite[p.~569]{AS3}
\[ch_g(\epsilon_M)=ch(\epsilon_M)=e^{c_1(L_M)/2}\prod^d_{a=1}(e^{p_a/2}-e^{-p_a/2})\]
($e^{c_1(L_M)/2}$ is the extra factor for the $spin^c$ structure, as formally $S_{spin^c}\cong L^{1/2}\otimes S_{Spin}$, also $g$ acts trivially on $\epsilon_M$). Working formally with the Chern roots in the universal
power series ring we obtain
\begin{equation}\label{localcom1}
 ch_g(S_N^+-S^-_N):=ch_g(\epsilon_N)=\frac{ch(\psi^n(\epsilon_M))}{ch(\epsilon_M)}   
=
e^{(n-1)c_1(L_M)/2}\prod_{a=1}^d
\frac{e^{np_a/2}-e^{-np_a/2}}{e^{p_a/2}-e^{-p_a/2}}
\end{equation}
By (\ref{N_of_symmetric})
\[N_{\mathbb{C}}=TM_{\mathbb{C}}\otimes(\rho_n-1)_{\mathbb{C}}=\bigoplus^{n-1}_{m=1} TM_{\mathbb{C}}\otimes\chi_m\]
\[\lambda_{-1}N_{\mathbb{C}}=\prod^{n-1}_{m=1} \lambda_{-1}(TM_{\mathbb{C}}\otimes\chi_m)\]
Taking equivariant Chern character, each summand contributes two factors $\prod^{d}_{a=1}(1-\zeta^me^{p_a})$ and $\prod^{d}_{a=1}(1-\zeta^me^{-p_a})$, where $\zeta^m=e^{2\pi im/n}$. So
\[ch_g(\lambda_{-1}N_{\mathbb{C}})=\prod^{d}_{a=1}\prod^{n-1}_{m=1} (1-\zeta^me^{p_a})(1-\zeta^me^{-p_a})\]
Using the cyclotomic identity, we get
\begin{equation}\label{localcom2}
    ch_g(\lambda_{-1}N_{\mathbb{C}})=\prod^{d}_{a=1}\frac{(1-e^{np_a})(1-e^{-np_a})}{(1-e^{p_a})(1-e^{-p_a})}=\prod^{d}_{a=1}(\frac{e^{np_a/2}-e^{-np_a/2}}{e^{p_a/2}-e^{-p_a/2}})^2
\end{equation}
Combining (\ref{localcom1}) and (\ref{localcom2}), we get
\[\frac{ch_g(S^{+}_N-S^{-}_N)}{ch_g(\lambda_{-1}(N_{\mathbb{C}}))}=e^{(n-1)c_1(L_M)/2}\prod^d_{a=1}  \frac{e^{p_a/2}-e^{-p_a/2}}{ (e^{np_a/2}-e^{-np_a/2})}\]
\end{proof}

\begin{lemma}\label{functoriallocalization} 
    For a continuous map $f:M\rightarrow X$  of compact CW-complexes, and $g=\sigma_n , G=\langle g\rangle$, the following diagram is commutative
\begin{equation}\label{diag3proof}
    \begin{tikzcd}
     K^{G,an}_0(M^n)_{g}  \arrow[r,  "\gamma^M_{g}"] \arrow[d,  "(f^n)_*"]& K^{G,an}_0(\Delta(M))_{g} \arrow[d,  "f_*"]   \\
K^{G,an}_0(X^n)_{g} \arrow[r,  "\gamma^X_{g}"] & K^{G,an}_0(\Delta(X))_{g}
     \end{tikzcd}
\end{equation}
\end{lemma}
\begin{proof}
Let  $\Delta(M)\xrightarrow{i_M}M^n, \Delta(X)\xrightarrow{i_X}X^n$. Recall 
$\gamma^M_g=(i_{M})_*^{-1}$ where $(i_{M})_*$ is defined by the Kasparov product with
$[i^*_M]\in KK^G(C(M^n),C(\Delta(M)))$. Similarly $\gamma^X_g=(i_{X})_*^{-1}$  where $(i_{X})_*$  is defined by  the Kasparov product with $[i^*_X]$.  Since $f^n\circ i_M=i_X \circ f$, the Kasparov product implies
\begin{equation}\label{kasppush}
    (f^n)_*\circ (i_{M})_*= (i_{X})_*\circ f_*
\end{equation}
which proves the commutativity of (\ref{diag3proof}).
\end{proof}

\textbf{Proof of Theorem \ref{analyticmoonen}}:

\begin{proof}
As mentioned in (\ref{f=g,bm=kk})
 \[\mathcal{F}=f_*(\mathcal{G})\]
where $\mathcal{G}=E\cap[\bD_M^{spin^c}] \in K_0^{an}(M)$, and $f:M\rightarrow X$ is a continuous map.  We have
\begin{equation}\label{CHcoommutef_*}
    Ch(P_n(\mathcal{F}),g)=f_*(Ch(P_n(\mathcal{G}),g))
\end{equation}
Indeed, $f^n:M^n\rightarrow X^n$ is automatically $S_n$-equivariant and
$ (X^n)^{g}=\Delta(X)\cong X, (M^n)^{g}=\Delta(M)\cong M$, and $f^n$ restricts on fixed loci to $f$.  Chern character and  homological Adams operation  $\Psi_n$ commute with proper pushforward. We also have
$(f_*\mathcal{G})^{\boxtimes n}=(f^n)_*(\mathcal{G}^{\boxtimes n})$. By Lemma \ref{functoriallocalization}, localization in this setting also commutes with $f$. Hence
    \begin{align*}
       f_*(Ch(P_n(\mathcal{G}),g))=f_*((Ch_g\circ\gamma^M_g)(\mathcal{G}^{\boxtimes n})_g)=&(Ch_g\circ\gamma^X_g)((f^n)_*(\mathcal{G}^{\boxtimes n})_g)=(Ch_g\circ\gamma^X_g)((f_*\mathcal{G})^{\boxtimes n}_g)\\
        =&Ch(\mathcal{F}^{\boxtimes n},g)=Ch(P_n(\mathcal{F}),g)
    \end{align*}
proving (\ref{CHcoommutef_*}). Then
\[
\begin{aligned}
Ch(P_n(\mathcal F),g)=f_*Ch(P_n(\mathcal G),g)=f_*\Psi_nCh(\mathcal G)=\Psi_nf_*Ch(\mathcal G)=\Psi_nCh(\mathcal F)
\end{aligned}
\]
where we used Theorem \ref{problemmoon} for $\mathcal{G}$.
\end{proof}

\section{Equivariant $K$-homology class of signature operator of $X^n$}

Albin-Banagl-Piazza \cite{ABP1} showed that the $K$-homology class of the signature operator in $K^{an}_0(X)$ is multiplicative, generalizing a result of Rosenberg-Weinberger \cite{RW2}. See also Remark \ref{anotherproofofcaretstionKhomology} for the equivariant case. 

    \begin{theorem}\label{signatureproduct} \cite[Corollary 5.18]{ABP1} Let $X$ be  an even dimensional oriented closed smoothly stratified Witt space. Then in $K^{an}_0(X^n)$ we have
    \begin{equation}\label{amp1externalsignture}
\left[D^{sign}_{X^n}\right]=\left[D^{sign}_{X}\right]^{\boxtimes n}
    \end{equation}
    where  $\boxtimes n$ denotes  the $n$-fold external Kasparov product.
\end{theorem}

By Theorem \ref{Sn-mather-X^n}, $X^n$ admits an $S_n$-Thom--Mather structure so by \cite{BLP} there is an equivariant $K$-homology class of the signature operator (assuming $\dim X$ is even) $$\left[D^{sign,S_n}_{X^n}\right]\in K_{\dim X}^{S_n,an}(X^n)\cong K_0^{S_n,an}(X^n)$$ which can be defined  by choosing an $S_n$-invariant wedge metric on $X^n$, e.g., the product metric.

On the other hand, $\left[D^{sign}_{X}\right]^{\boxtimes n}$ has a natural $S_n$-equivariant structure as we saw in the construction of the power operation, where $S_n$ acts by the graded permutation on factors of the external Kasparov product, which defines a class  $P_n(\left[D^{sign}_{X}\right]) \in K_0^{S_n,an}(X^n)$.

As mentioned in the introduction, in general, these two equivariant structures are \textit{not} the same. To compare $\left[D^{sign,S_n}_{X^n}\right]$  and $P_n(\left[D^{sign}_{X}\right])$ we record their representative by equivariant unbounded Kasparov cycles. 

We choose the metric $g$ on $X_{reg}$ to be an \textit{adapted rigid iterated edge metric} in the sense of \cite[Prop.~5.4]{ALMP}. In this case the signature operator $D$ initially defined on $\Omega^*_{c}(X_{reg})$ becomes essentially self-adjoint \cite[Thm.~1.1]{ALMP}, hence the self-adjoint extension is unique. Note that the $K$-homology class of signature operator $\left[D_X^{sign}\right]$ does not depend on the type of wedge metric \cite{ALMP}.  If
 \[\left[D_X^{sign}\right]=[(\mathcal{H}_X,
\tau_X,\rho_{X},
D_X)]\]
where $\mathcal{H}_X
=L^2(X_{reg},\Omega^*, g)$, then by (\ref{poweropanal}) we have
\begin{equation}\label{powersignkcycle}
    P_n(\left[D^{sign}_{X}\right])=\left[(\mathcal{H}_X^{\widehat{\otimes}n},
\tau_X^{\otimes n},\rho^{\widehat{\otimes}n}_{X},
D_{X,\tau}^{(n)},U^{n}_{\tau})\right]
\end{equation}

We equip $X^n$ with the product metric $g^n$. Let

\[\left[D^{sign,S_n}_{X^n}\right]=[(\mathcal{H}_{X^n},\tau_{X^n},\rho_{X^n},D_{X^n},  \widetilde{U})]\]
where    $D_{X^n}$ denotes the  signature operator of $X^n$ and
\[\mathcal{H}_{X^n}:=L^2(X^n_{reg},\Omega^*, g^n)\]
is the underlying Hilbert space, and 
\[\widetilde{U}:S_n\rightarrow \mathcal{U}(\mathcal{H}_{X^n})\]
is the corresponding  unitary representation where  by the definition of equivariant  $K$-homology class (\ref{equivrepdef}) for $\sigma\in S_n$  it is defined by 
\begin{equation}
    \widetilde{U}_{\sigma} w=(\sigma^{-1})^*w
\end{equation}
where $\sigma$ acts on $X^n$ by 
 \[
\sigma(x_1,\ldots,x_n)
=
(x_{\sigma^{-1}(1)},\ldots,x_{\sigma^{-1}(n)})
\]
we can describe the action 
$$  \widetilde{U}_{\sigma}=(\sigma^{-1})^*:\mathcal{H}_{X^n} \rightarrow \mathcal{H}_{X^n}$$
by using the identification (a unitary isomorphism)
\begin{equation}\label{phiidentif}
    \phi: \mathcal{H}_X^{\widehat{\otimes}n}=\mathcal{H}_X\widehat{\otimes}
    \cdots\widehat{\otimes}\mathcal{H}_X \xrightarrow{\cong}\mathcal{H}_{X^n}
\end{equation}
which on smooth compactly supported forms $\alpha_i\in \Omega^{p_i}_{c}(X_{reg})$ is defined by  
\[(\alpha_1 \widehat{\otimes} \ldots \widehat{\otimes}\alpha_n)\mapsto \alpha_1\boxtimes \ldots \boxtimes \alpha_n
\]
\[\alpha_1\boxtimes \ldots \boxtimes \alpha_n:=\pi^*_1(\alpha_1)\wedge \ldots \wedge \pi^*_n(\alpha_n)\in \Omega_c^{\sum^n_{i=1} {p_i}}(X^n_{reg})\]
where $\pi_i$ is the projection to the $i$-th coordinate. The map is isometric on elementary tensors because the metric used for $\mathcal{H}_{X^n}$ is the product metric $g^n$. Then since finite sums of exterior products of compactly supported smooth forms are dense in $L^2\Omega^*(X^n_{\mathrm{reg}}, g^n)$, $\phi$  extends to a unitary isomorphism on $ \mathcal{H}_X^{\widehat{\otimes}n}$.


\begin{lemma}
    Let $X$ be an even dimensional oriented closed smoothly stratified Witt space. Equip $X$ with a wedge metric $g$, and $X^n$ with the product metric $g^n$ . Then
    \begin{enumerate}
        \item $\phi:\mathcal{H}_X^{\widehat{\otimes}n}\xrightarrow{\cong}\mathcal{H}_{X^n} $ is a unitary isomorphism of \(\mathbb Z_2\)-graded Hilbert spaces where the gradings are $\tau^{\otimes n}_X$ and $\tau_{X^n}$ respectively, i.e., 
        \[\phi^{-1} \tau_{X^n}\phi=\tau^{\otimes n}_X\]
        \item  for $\sigma\in S_n$ \[\phi^{-1} \widetilde{U}_\sigma \phi=U^n_{\epsilon,\sigma}\]
         recall the notation \begin{equation}\label{geosign}
          U^n_{\epsilon,\sigma}(\xi_1\widehat{\otimes}\cdots
    \widehat{\otimes}\xi_n)
    =
    (-1)^{\epsilon(\sigma;\xi)}
    \xi_{\sigma^{-1}(1)}\widehat{\otimes}
    \cdots\widehat{\otimes}
    \xi_{\sigma^{-1}(n)}
    \end{equation}
    where
    \[
    \epsilon(\sigma;\xi)
    =
    \sum_{\substack{i<j\\
    \sigma(i)>\sigma(j)}}
    |\xi_i|_{\epsilon}|\xi_j|_{\epsilon}
    \quad(\mathrm{mod}\ 2)\]
   where $|\alpha|_{\epsilon}=0$ if the degree (as a differential form) of $\alpha$ is even, otherwise it is $1$.
    \end{enumerate}
\end{lemma}
\begin{proof}
The first assertion is basic calculation from the definition (\ref{siggradingdef}), we mention the relation for  Hodge $\ast$-operators in the case of $n=2$ for simplicity. For a $p$-form $\alpha \in \Omega^{p}_{c}(X_{reg})$ and  a $q$-form $\beta \in \Omega^{q}_{c}(X_{reg})$ 
\[\ast_{X\times X}(\alpha \boxtimes\beta)=(-1)^{(2m-p)q}\left(\ast_X(\alpha) \boxtimes \ast_X(\beta) \right)\]
where  $\dim X_{reg}=2m$.
From this one obtains 
\[\tau_{X^2}(\phi(\alpha {\otimes} \beta))=\tau_{X\times X}(\alpha \boxtimes\beta)=\tau_X(\alpha) \boxtimes \tau_X(\beta)=\phi((\tau_X\otimes\tau_X )(\alpha {\otimes} \beta))\]
To prove the second assertion, it suffices to check this for a $\sigma$ being a transposition $\sigma=(1,2)$, then the parity-of-degree grading is the sign that comes from graded commutativity of the wedge product;
\begin{align*}
    \widetilde{U}_\sigma (\phi(\alpha {\otimes} \beta))=\widetilde{U}_\sigma(\alpha \boxtimes\beta)=(\sigma^{-1})^*(\pi^*_1(\alpha)\wedge \pi^*_2(\beta))=&\pi^{*}_2(\alpha)\wedge \pi_1^*(\beta)=(-1)^{pq}\pi^{*}_1(\beta)\wedge \pi_2^*(\alpha)\\=&
    (-1)^{pq}\beta \boxtimes \alpha=\phi((-1)^{|\alpha|_\epsilon|\beta|_\epsilon} \beta\otimes \alpha)\\=&\phi U^2_{\epsilon,\sigma}(\alpha {\otimes}\beta)
\end{align*}
\end{proof}

\begin{lemma}
    Let $X$ be an even dimensional oriented closed smoothly stratified Witt space. Equip $X_{reg}$ with an {adapted rigid iterated edge metric} $g$ in the sense of \cite{ALMP}, and $X^n$ with the product metric $g^n$. Then
   \begin{equation}\label{intertwindeD}
       \phi^{-1} D_{X^n} \phi=D^{(n)}_{X,\epsilon}
   \end{equation}
        \begin{equation}\label{Dxepsilonn}
    D_{X,\epsilon}^{(n)}:=\sum_{j=1}^n
\epsilon^{\otimes(j-1)}
\otimes D_X \otimes 1^{\otimes(n-j)}
        \end{equation}
        \end{lemma}
\begin{proof}
For simplicity of notation  again assume $n=2$, the $n$-fold result follows by iteration.
First, we show the equality on  
 compactly supported smooth homogeneous forms $\alpha \in \Omega^{p}_{c}(X_{reg}), \beta \in \Omega^{q}_{c}(X_{reg})$.
  Since the metric is a product metric, we have
\[d_{X\times X}(\alpha\boxtimes\beta)=(d_X\alpha)\boxtimes\beta
+(-1)^p\alpha\boxtimes(d_X\beta)\]
\[\delta_{X\times X}(\alpha\boxtimes\beta)
=
(\delta_X\alpha)\boxtimes\beta
+
(-1)^p\alpha\boxtimes(\delta_X\beta)\]

\begin{align*}
     D_{X^2} \phi(\alpha {\otimes} \beta)=D_{X\times X}(\alpha\boxtimes\beta)=&(d_{X^2}+\delta_{X^2})(\alpha\boxtimes\beta)=(D_X\alpha)\boxtimes\beta
+
(-1)^p\alpha\boxtimes(D_X\beta)\\
=&\phi\left((D_X\otimes1 + \epsilon \otimes D_X)(\alpha {\otimes} \beta)\right)
=\phi D^{(2)}_{X,\epsilon}(\alpha \otimes \beta)
\end{align*}
therefore
\begin{equation}\label{eqinsidethm}
    \phi D^{(2)}_{X,\epsilon}  \phi^{-1}=D_{X^2} \qquad  on \qquad \phi(\Omega^*_{c}(X_{reg})\otimes \Omega^*_{c}(X_{reg})) \subset \Omega^*_{c}(X^2_{reg})
\end{equation}
note that $\phi(\Omega^*_{c}(X_{reg})\otimes \Omega^*_{c}(X_{reg}))$ 
 is not equal to $\Omega^*_{c}(X^2_{reg})$ in general.

  By the assumption on the metric  and \cite[Thm~1.1]{ALMP}, $D_X$ is an   essential self adjoint operator initially defined on
$\Omega_c^*(X_{reg})$.  Then by \cite[Prop.~3.6 and Cor.~3.7]{SV} the operator  $D^{(2)}_{X,\epsilon}$ in (\ref{Dxepsilonn})  is also essential self adjoint on the algebraic  tensor product  $\Omega^*_{c}(X_{reg})\otimes \Omega^*_{c}(X_{reg})$ (as $D_X$ and $\epsilon$ anticommute and the commutator hypotheses of loc. cit. are trivially satisfied). Since $\phi$ is unitary, $ \phi D^{(2)}_{X,\epsilon}  \phi^{-1}$ is  essentially self-adjoint on $\phi(\Omega^*_{c}(X_{reg})\otimes \Omega^*_{c}(X_{reg}))$. Since $D_{X^2}$ is a symmetric operator initially defined on $\Omega^*_{c}(X^2_{reg})$, (\ref{eqinsidethm}) implies that $D_{X^2}$ is also essential self adjoint, and taking closures implies (\ref{eqinsidethm}) holds on the domain of the closure.
\end{proof}

Recall  that by definition of $KK$-homology, the class of    two $K$-cycles which are unitarily equivalent are the same, i.e.,  if there is a unitary between their Hilbert spaces intertwining the grading, representations, and the operators. The previous two lemmas about the  unitary map $\phi$ imply

\begin{corollary}
    Let $X$ be an even dimensional oriented closed smoothly stratified Witt space. Then 
    \begin{equation}\label{cartesiansignkcycle}
    \left[D^{sign,S_n}_{X^n}\right]=\left[ (\mathcal{H}_{X^n},\tau_{X^n},\rho_{X^n},D_{X^n},  \widetilde{U})   \right]  =\left[(\mathcal{H}_X^{\widehat{\otimes}n},
\tau_X^{\otimes n},\rho^{\widehat{\otimes}n}_{X},
D_{X,\epsilon}^{(n)},U^n_{\epsilon})\right]
\end{equation}
\end{corollary}

\begin{remark}\label{anotherproofofcaretstionKhomology}
By a similar argument, one can show an equivariant version of \cite[Corollary 5.18]{ABP1}; namely
\begin{equation}\label{equivaintabp1result}
    [D^{sign, G_1\times G_2}_{X_1\times X_2}]=[D^{sign, G_1}_{X_1}] \boxtimes [D^{sign, G_2}_{X_2}]
\end{equation}
for even dimensional compact $G_1$-Witt space $X_1$ and $G_2$-Witt space $X_2$.  Indeed, we can represent for $a\in \{1,2\}$
\[[D^{sign, G_a}_{X_a}]=\left[(\mathcal{H}_{X_a},\tau_{X_a},\rho_{X_a},D_{X_a},{U}_{G_a})\right]\]
then
\begin{equation}\label{cycle1}
   [D^{sign, G_1}_{X_1}] \boxtimes [D^{sign, G_2}_{X_2}]=\left[(\mathcal{H}_{X_1}\widehat{\otimes} \mathcal{H}_{X_2},\tau_{X_1}\otimes \tau_{X_2},\rho_{X_1} \widehat{\otimes}\rho_{X_2},D_{X_1}\otimes 1+ \tau_{X_1}\otimes D_{X_2} ,{U}_{G_1}\times {U}_{G_2})\right] 
\end{equation}
\begin{equation}\label{cycle2}
    [D^{sign, G_1\times G_2}_{X_1\times X_2}]=\left[(\mathcal{H}_{X_1\times X_2},\tau_{X_1\times X_2},\rho_{X_1\times X_2},D_{X_1\times X_2},  {U}_{G_1\times G_2})    \right]
\end{equation}
Then the map 
\[\phi:\mathcal{H}_{X_1}\widehat{\otimes} \mathcal{H}_{X_2} \rightarrow \mathcal{H}_{X_1\times X_2}\]
is a $G_1\times G_2$-\textit{equivariant} unitary isomorphism of  \(\mathbb Z_2\)-graded Hilbert spaces where
\[
\phi^{-1}\tau_{X_1\times X_2}\phi=\tau_{X_1}\otimes\tau_{X_2}, \quad 
\phi^{-1}U_{(g_1,g_2)}\phi=U_{g_1}\otimes U_{g_2}
\]
while
\[\phi^{-1}D_{X_1\times X_2}\phi=D_{X_1}\otimes1+\epsilon_{X_1}\otimes D_{X_2}\]
Then  there is a unitary isomorphism  $W:\mathcal{H}_{X_1}\widehat{\otimes} \mathcal{H}_{X_2} \rightarrow \mathcal{H}_{X_1}\widehat{\otimes} \mathcal{H}_{X_2}$ defined by

\[
W:=
\begin{cases}
id\otimes id,&\text{on }\mathcal H_{X_1,0}\otimes\mathcal H_{X_2}\\
id\otimes\tau_{X_2},&\text{on }\mathcal H_{X_1,1}\otimes\mathcal H_{X_2}
\end{cases}
\]
where $\mathcal H_{X_1}=\mathcal H_{X_1,0}\oplus\mathcal H_{X_1,1}$ is the splitting with respect to $\kappa_1=\epsilon_{X_1}\tau_{X_1}$; $$\kappa_1|_{\mathcal H_{X_1,0}}=1, \qquad
\kappa_1|_{\mathcal H_{X_1,1}}=-1$$ Then
\[W(D_{X_1}\otimes 1+ \tau_{X_1}\otimes D_{X_2})W^{-1}=D_{X_1}\otimes1+\epsilon_{X_1}\otimes D_{X_2},\quad W(\tau_{X_1}\otimes \tau_{X_2})W^{-1}=\tau_{X_1}\otimes \tau_{X_2}\]
therefore $\phi W$ intertwines the cycles (\ref{cycle1}) and (\ref{cycle2}).
\end{remark}

\subsection{Proof of Theorem \ref{introtheoreminductform}}
To simplify the notation we drop the subscript $X$ from $D_X, \cH_X,\tau_X,\rho_X$. 
  \[\left[D_X^{sign}\right]=[(\mathcal{H},
\tau,\rho,
D)]\]
\begin{proof}
 We have in $K_0^{an}(X^n)$
\begin{align}
    \left[D^{sign}_{X^n}\right]=[D^{sign}_X]^{\boxtimes n}=([\tilde{S}_0] -[\tilde{S}_1])^{\boxtimes n}=&\sum_{\mathfrak{a} \in \{0,1\}^n} (-1)^{|\mathfrak{a}|}[\tilde{S}_{a_1}]\boxtimes \ldots \boxtimes [\tilde{S}_{a_n}]\\ \label{nonequivnatequalitydx^n}
    =&\sum^n_{k=0}(-1)^k \left( \sum_{|\mathfrak{a}|=k}  [\tilde{S}_{a_1}]\boxtimes \ldots \boxtimes [\tilde{S}_{a_n}] \right)
\end{align}
where $\mathfrak{a}=(a_1,\ldots,a_n)\in \{0,1\}^n, |\mathfrak{a}|=\sum a_i$. Let
\[\cH^{\widehat{\otimes}\mathfrak{a}}:=\cH_{a_1}\widehat{\otimes} \ldots \widehat{\otimes} \cH_{a_n},\quad \tau^{{\otimes}\mathfrak{a}}:=\tau_{a_1}\otimes\ldots\otimes \tau_{a_n}, \quad \rho^{\widehat{\otimes}\mathfrak{a}}:=\rho_{a_1}\widehat{\otimes} \ldots \widehat{\otimes} \rho_{a_n}\]
\begin{equation}\label{D^a_eps}
D_\epsilon^{\mathfrak a}:=
\sum_{j=1}^n
\epsilon_{a_1}\otimes\cdots\otimes\epsilon_{a_{j-1}}
\otimes D_{a_j}\otimes
1^{\otimes(n-j)}
\end{equation}
 e.g., for $n=2, D_{\epsilon}^{(a_1,a_2)}=D_{a_1}\otimes1+\epsilon_{a_1} \otimes D_{a_2} $. Define the equivariant $K$-cycle
\begin{equation}\label{defD_k}
    \cD_k=(\bigoplus_{{|\mathfrak{a}|=k}} \cH^{\widehat{\otimes}\mathfrak{a}},\bigoplus_{{|\mathfrak{a}|=k}} \tau^{{\otimes}\mathfrak{a}},\bigoplus_{{|\mathfrak{a}|=k}} \rho^{\widehat{\otimes}\mathfrak{a}},\bigoplus_{{|\mathfrak{a}|=k}}D_{\epsilon}^{\mathfrak{a}}, U_{k,\epsilon})
\end{equation}
\[U_{k,\epsilon}:S_n\rightarrow \mathcal{U}(\bigoplus_{{|\mathfrak{a}|=k}} \cH^{\widehat{\otimes}\mathfrak{a}})\]
where for $\sigma\in S_n,\xi_i\in \cH_{a_i}, (\xi_1\widehat{\otimes}\cdots
    \widehat{\otimes}\xi_n)\in \cH^{\widehat{\otimes}\mathfrak{a}}$ we have
    \begin{equation}\label{equivaction}
        (U_{k,\epsilon})_{\sigma}(\xi_1\widehat{\otimes}\cdots
    \widehat{\otimes}\xi_n)=(-1)^{\epsilon(\sigma;\xi)}
    \xi_{\sigma^{-1}(1)}\widehat{\otimes}
    \cdots\widehat{\otimes}
    \xi_{\sigma^{-1}(n)}\in \cH^{\widehat{\otimes}\sigma\mathfrak{ a}}
    \end{equation}
$U_{k,\epsilon}$ commutes with the operator $\bigoplus_{{|\mathfrak{a}|=k}}D_{\epsilon}^{\mathfrak{a}}$ and conditions of Definition \ref{equivkkdef} are satisfied, therefore (\ref{defD_k}) defines a class $[\cD_k]\in K^{S_n,an}_0(X^n)$. Then we have the equality in $K^{S_n,an}_0(X^n)$
\[\left[D^{sign,S_n}_{X^n}\right]=\left[(\mathcal{H}^{\widehat{\otimes}n},
\tau^{\otimes n},\rho^{\widehat{\otimes}n}_{},
D_{\epsilon}^{(n)},U^n_{\epsilon})\right]=\sum^n_{k=0} [\cD_k]\]
where the first equality is (\ref{cartesiansignkcycle}) and the second equality comes from the fact that
\[U^n_{\epsilon}=\bigoplus^n_{k=0}U_{k,\epsilon}\]
note that for each $|\mathfrak{a}|=k$, we have $\tau^{{\otimes}\mathfrak{a}}=(-1)^k \epsilon^{{\otimes}\mathfrak{a}}$. Therefore
 as a class in $K^{S_n,an}_0(X^n)$
\[[\cD_k]=(-1)^k\left[ (\bigoplus_{{|\mathfrak{a}|=k}} \cH^{\widehat{\otimes}\mathfrak{a}},\bigoplus_{{|\mathfrak{a}|=k}} \epsilon^{{\otimes}\mathfrak{a}},\bigoplus_{{|\mathfrak{a}|=k}} \rho^{\widehat{\otimes}\mathfrak{a}},\bigoplus_{{|\mathfrak{a}|=k}}D_{\epsilon}^{\mathfrak{a}}, U_{k,\epsilon})\right]\]
we claim that
\begin{equation}\label{induecdk}
  \left[ (\bigoplus_{{|\mathfrak{a}|=k}} \cH^{\widehat{\otimes}\mathfrak{a}},\bigoplus_{{|\mathfrak{a}|=k}} \epsilon^{{\otimes}\mathfrak{a}},\bigoplus_{{|\mathfrak{a}|=k}} \rho^{\widehat{\otimes}\mathfrak{a}},\bigoplus_{{|\mathfrak{a}|=k}}D_{\epsilon}^{\mathfrak{a}}, U_{k,\epsilon})\right]=Ind_{S_{n-k}\times S_k}^{S_n}(P_{n-k}([\widetilde S_{0}])\boxtimes P_{k}[\widetilde S_{1}])
\end{equation}
see subsection \ref{inductkcycles} for the notion of induction on the level of $K^{an}_0$. Fix the element
\[
\mathfrak a_0=
(\underbrace{0,\ldots,0}_{n-k},
\underbrace{1,\ldots,1}_{k})
\]
\[
\mathcal H^{\widehat\otimes\mathfrak a_0}
=
\mathcal H_0^{\widehat\otimes(n-k)}
\widehat\otimes
\mathcal H_1^{\widehat\otimes k}
\]
recall that for $a\in \{0,1\}$, $\widetilde S_a=(\mathcal H_a,\epsilon_a,\rho_a,D_a)$. We have by definition of power operation and Kasparov external product
\begin{equation}\label{poweropn-k,k}
    P_{n-k}([\widetilde S_0])
\boxtimes P_k([\widetilde S_1])=[\left(
\mathcal H^{\widehat\otimes\mathfrak a_0},\epsilon^{\otimes\mathfrak a_0},
\rho^{\widehat\otimes\mathfrak a_0},
D_\epsilon^{\mathfrak a_0},
U^{n-k}_{\epsilon}\times U^{k}_{\epsilon}
\right)] 
\end{equation}
where $U^{n-k}_{\epsilon}\times U^{k}_{\epsilon}$ is the graded permutation action of $S_{n-k}\times S_k$ constructed from $$U^{n-k}_{\epsilon}:S_{n-k}\rightarrow \mathcal{U}(\cH_0^{\widehat\otimes(n-k)}), \quad
 U^{k}_{\epsilon}:S_{k}\rightarrow \mathcal{U}(\cH_1^{\widehat\otimes(k)})$$
 The stabilizer of the multi-index \(\mathfrak a_0\) is \(S_{n-k}\times S_k\), and its orbit indexes the summands, i.e., there is a bijection 
\[
S_n/(S_{n-k}\times S_k)
\longrightarrow
\{\mathfrak a\in\{0,1\}^n:|\mathfrak a|=k\},
\qquad
\sigma(S_{n-k}\times S_k)\longmapsto
\sigma\mathfrak a_0
\]
Then from the orbit-stabilizer explanation  in subsection \ref{inductkcycles} and (\ref{indcutfomrulakcycles}) it is clear that by applying $Ind_{S_{n-k}\times S_k}^{S_n}$ to  the right hand side of (\ref{poweropn-k,k}) we get the left hand side of (\ref{induecdk}), as the gradings and equivariant structures agree as just explained. Therefore by combining everything we get
\begin{equation}
    \left[D^{sign,S_n}_{X^n}\right]=\sum^n_{k=0}(-1)^kInd_{S_{n-k}\times S_k}^{S_n}(P_{n-k}([\tilde{S}_{0}])\boxtimes P_{k}([\tilde{S}_{1}]))
\end{equation}
(with the convention $P_0(x)=1$).
Note that by applying the forgetful functor we get (\ref{nonequivnatequalitydx^n}).

\end{proof}

\subsection{Equivariant $L$-classes of $X^n$}
As an application we can compute the localized Chern character $Ch( [D^{sign,S_n}_{X^n}],\sigma_n)$ of an $n$-cycle $\sigma_n\in S_n$.

\begin{theorem} \label{vanishingchind}  Following the notation used in Theorem \ref{introtheoreminductform} we have:
\begin{enumerate} 
\item for $i\in \{0,1\}$ we have
\[Ch( P_{n}([\tilde{S}_{i}]),\sigma_n)=\Psi_nCh([\tilde{S}_{i}])\]
    \item for $0<k<n$ we have
    \[Ch( Ind_{S_{n-k}\times S_k}^{S_n}(P_{n-k}([\tilde{S}_{0}])\boxtimes P_{k}([\tilde{S}_{1}])),\sigma_n)=0\]
    
\end{enumerate}
\end{theorem}
\begin{proof}
(1) follows from  Theorem \ref{analyticmoonen}. (2) follows from applying the representation theoretic Lemma \ref{mackylem} to the definition of localized Chern character $Ch(-,\sigma_n)$,  Definition \ref{localizedchernclass}. Indeed, to compute $Ch(-,\sigma_n)$ we start by applying the maps
\[K_0^{S_n,an}(X^n)\xrightarrow{res^{S_n}_{\langle \sigma_n \rangle}} K_0^{\langle \sigma_n \rangle,an}(X^n)\rightarrow K_0^{\langle \sigma_n \rangle,an}(X^n)_{{\sigma_n}}\]
If    $0<k<n$, no conjugate of an $n$-cycle $\sigma_n$ belongs to $S_{n-k}\times S_k$, as $\sigma_n$ does not preserve any nonempty proper subset of $\{1,\ldots,n\}$.
\end{proof}

\begin{lemma}\label{mackylem} Let $G$ be a finite group acting on a locally compact space $X$.
    Let $H\subset G$ be a subgroup. Let $g\in G$ such that no conjugate of $g$ lies in $H$. Then for any $x\in K_0^{an,H}(X)$
    \begin{equation}\label{eqlemmavanishignlocal}
        \left( res^G_{\langle g\rangle}(Ind^G_H(x)) \right)_g=0
    \end{equation}
\end{lemma}
\begin{proof}
Recall $R(G)$ denotes the complex representation ring of $G$. By assumption no conjugate of $g$ lies in $H$, therefore there is an element $a_g\in R(G)$ whose character at $g$ is non-zero, while  $\chi_{a_g}(h) =0$, when $h\in H$ (apply \cite[Prop.~3.7(iv)]{Se} to the prime ideal $\frak{m}_g$ consisting of virtual representations whose character vanishes at $g$). Hence $res^G_H(a_g)=0$  . By the projection formula (\ref{projectionfomrk})
\[a_g.Ind^G_H(x)=Ind^G_H(res^G_H(a_g).x)=0\]
applying $res^G_{\langle g\rangle}$ we get
\[res^G_{\langle g\rangle}(a_g).res^G_{\langle g\rangle}(Ind^G_H(x))=0\]
by construction $res^G_{\langle g\rangle}(a_g)$ is invertible in the representation ring localized at $g$, therefore we get (\ref{eqlemmavanishignlocal}) after localization.
\end{proof}

Theorems \ref{introtheoreminductform} and \ref{vanishingchind} imply Theorem \ref{ch_sign_x^n-zagier}, which in turn implies Theorem \ref{equianvtLclasszagiersingular}:

\begin{corollary}\label{ch_sgiamr_x^r}
Assume \(X\) is connected, then
\[
Ch(\left[D^{sign,S_n}_{X^n}\right],\sigma_n)
=
\begin{cases}
\Psi_nCh([\tilde{S}_{0}]+[\tilde{S}_{1}])
=\Psi_nCh([D^{dR}_X])=\Psi_nI\chi(X)=I\chi(X), & \text{if } n \text{ is even},\\[4pt]
\Psi_nCh([\tilde{S}_{0}]-[\tilde{S}_{1}])
=\Psi_nCh([D^{sign}_X])=\Psi_nL^{an}_*(X), & \text{if } n \text{ is odd}.
\end{cases}
\]
\end{corollary}

\section{The characteristic class Hodge index conjecture for symmetric products $X^{(n)}$}

This section is devoted to the proof of Theorems \ref{symLclassthm} and \ref{bsyconjforsymmetric}. The canonical identification $X^n\times X^m\rightarrow X^{n+m}$ induces
the Pontryagin product on the symmetric products  \cite[Def. 1.1]{CMSSY} 
\[X^{(n)}\times X^{(m)} \xrightarrow{} X^{(n+m)}\]
i.e., by "concatenating"  unordered tuples of points  and by the embedding $S_n\times S_m \hookrightarrow S_{n+m}$.
 Let $F$ be the Borel-Moore homology functor  defined on complex
quasi-projective varieties. By functoriality this also induces
\[F(X^{(n)})\times F(X^{(m)})\xrightarrow{\boxtimes}F(X^{(n)}\times X^{(m)}) \xrightarrow{(-)_*} F(X^{(n+m)})\]
Given an element $\sigma\in S_n$, let the cycle decomposition be $\sigma=\prod_r\prod^{k_r}_{j=1}\sigma_{r,j}$ where  $\sigma_{r,j}$ is an $r$-cycle. For a cycle $A$, denote by $|A|$ its length. We have \cite[\S~3.1]{CMSSY}
\[(X^n)^{\sigma}\cong \prod_{r,j} (X^{|\sigma_{r,j}|})^{\sigma_{r,j}}\cong \prod_{r}\Delta_r(X)^{k_r}\cong X^{k_1+k_2+\ldots}\]
 The inclusion $X^{|\sigma_{r,j}|}\hookrightarrow X^n$ is given by $X_{l_1}\times \ldots X_{l_r}$, where $\sigma_{r,j}=(l_1,\ldots,l_r)$, with $X_{l_i}$ on the $l_i$-th place in $X^n$. $\Delta_r(X)$ is the diagonal in $X^r$.  The projection $\pi^{\sigma}:(X^n)^{\sigma}\rightarrow X^n \rightarrow X^{(n)}$ is the product of
\[\pi^{\sigma_{r,j}}:X\cong (X^{|\sigma_{r,j}|})^{\sigma_{r,j}} \hookrightarrow X^{|\sigma_{r,j}|} \rightarrow X^{(|\sigma_{r,j}|)}\]
with $\prod_{r,j} X^{(|\sigma_{r,j}|)}\rightarrow X^{(n)}$
induced by the Pontryagin product.

\begin{theorem}\label{symmetricfactorchern}
Let $X$ be a compact oriented even-dimensional
smoothly stratified Witt space, and let
$\sigma\in S_n$ with the cycle decomposition $\sigma=\prod_r\prod^{k_r}_{j=1}\sigma_{r,j}$ where  $\sigma_{r,j}$ is an $r$-cycle. Then
    \[Ch(\left[D^{sign,S_n}_{X^n}\right],\sigma)=\prod_{r,j}  Ch([D^{sign,S_r}_{X^r}],\sigma_{r,j}) \in  \bigotimes_{r,j} H_*((X^{r})^{\sigma_{r,j}},\mathbb{C})\cong H_*((X^n)^{\sigma},\mathbb{C})\]
  where  $\prod_{r,j}$ denotes the external
product of homology classes.
 \end{theorem}
\begin{proof}
The proof follows from the following facts. Let $G':=\langle \sigma\rangle ,H:=\prod_r \prod^{k_r}_{j=1} \langle  \sigma_{r,j}\rangle$. To calculate $Ch(\left[D^{sign,S_n}_{X^n}\right],\sigma)$, by definition, we first restrict to $G'$ from $S_n$ and then apply the localization map and then the equivariant Chern character:
  \begin{equation}\label{firstrow}
     K_0^{S_n,an}(X^n)\xrightarrow{res^{S_n}_{G'}} K_0^{G',an}(X^n)\rightarrow K_0^{G',an}(X^n)_{{\sigma}}\xrightarrow[]{\gamma^{G'}_{\sigma}}K_0^{G',an}((X^n)^\sigma)_{\sigma}=K^{an}_0((X^n)^\sigma) \otimes R(G')_{\sigma}
  \end{equation}
  \[ \xrightarrow{Ch^{(X^n)^{\sigma}}\otimes \chi_\sigma}H_*((X^n)^\sigma,\mathbb{C}) \]
    We obtain the same class if we restrict to $H$ instead of $G'\subset H$. Indeed, we have
    \begin{equation}\label{secondrow}
        K_0^{S_n,an}(X^n)\xrightarrow{res^{S_n}_{H}} K_0^{H,an}(X^n)\rightarrow K_0^{H,an}(X^n)_{{\sigma}}\xrightarrow[]{\gamma^{H}_{\sigma}}K_0^{H,an}((X^n)^\sigma)_{\sigma}=K^{an}_0((X^n)^\sigma) \otimes R(H)_{\sigma}
    \end{equation}
\[\xrightarrow{Ch^{(X^n)^{\sigma}}\otimes \chi_\sigma}H_*((X^n)^\sigma,\mathbb{C})\]
where we used the fact that $H$ acts trivially on $(X^n)^{\sigma}$. Now, $\chi_\sigma:R(H)_{\sigma} \rightarrow \mathbb{C}$ factors through $R(H)_{\sigma}\xrightarrow{\operatorname{res}^{H}_{G'}} R(G')_{\sigma}\xrightarrow{\chi_\sigma}\mathbb{C}$, and there are equivariant vertical natural restriction maps from   (\ref{secondrow}) to (\ref{firstrow}), then we use the compatibility  $
\operatorname{res}^{H}_{G'}
\circ\gamma_\sigma^H
=
\gamma_\sigma^{G'}
\circ\operatorname{res}^{H}_{G'}
$.
We  have 
\[res^{S_n}_H \left[D^{sign,S_n}_{X^n}\right]=\boxtimes_{r,j} res_{\langle  \sigma_{r,j}\rangle}^{S_r} [D^{sign,S_r}_{X^r}]\]
which follows from the equivariant multiplicativity  of signature operators (\ref{equivaintabp1result}).
  The localization map is multiplicative with respect to the external product and product of groups (by the definition of localization and functoriality of Kasparov product)
\[\gamma_{\sigma}\left(\boxtimes_{r,j} res_{\langle  \sigma_{r,j}\rangle}^{S_r} [D^{sign,S_r}_{X^r}]\right)_{\sigma}=\boxtimes_{r,j} \gamma^{\langle \sigma_{r,j} \rangle}_{{\sigma_{r,j}}}(res_{\langle  \sigma_{r,j}\rangle}^{S_r} [D^{sign,S_r}_{X^r}])_{\sigma_{r,j}}\]
where $\gamma^{\langle \sigma_{r,j} \rangle}_{{\sigma_{r,j}}}:K_0^{\langle  \sigma_{r,j}\rangle,an}(X^{|\sigma_{r,j}|})_{\sigma_{r,j}} \rightarrow K_0^{\langle  \sigma_{r,j}\rangle,an}((X^{|\sigma_{r,j}|})^{\sigma_{r,j}})_{\sigma_{r,j}}$.
 Both $Ch$ and $\chi_{\sigma}$ satisfy multiplicativity under  K\"unneth isomorphism
\[H_*((X^n)^{\sigma},\mathbb{C}) \cong \bigotimes_{r,j} H_*((X^{|\sigma_{r,j}|})^{\sigma_{r,j}},\mathbb{C})\]
By combining everything we get the result.
\end{proof}

\subsection{Proof of Theorem \ref{symLclassthm}}
\begin{proof}
We follow \cite{CMSSY}. 
\begin{align*}
    L^{an}_*(X^{(n)}) \xlongequal{\text{Corollary~\ref{averageLclass}}} &
    \frac{1}{|S_n|}\sum_{\sigma\in S_n} \pi^{\sigma}_* L^{an}_*(X^n,\sigma)=\frac{1}{|S_n|}\sum_{\sigma\in S_n} \pi^{\sigma}_*(Ch(\left[D^{sign,S_n}_{X^n}\right],\sigma))\\
    \xlongequal{\text{Theorem~\ref{symmetricfactorchern}}} & 
    \frac{1}{|S_n|}\sum_{\sigma\in S_n}  \prod_{r,j}(\pi^{\sigma_{r,j}}_* Ch([D^{sign,S_r}_{X^r}],\sigma_{r,j}))\\
\xlongequal{\text{Corollary~\ref{ch_sgiamr_x^r}}} & 
    \frac{1}{n!}\sum_{\sigma\in S_n}  \prod_{r,j} (\pi^{\sigma_{r,j}}_*\Psi_r C_{(-1)^r}(X)) \in H_*(X^{(n)},\mathbb{C})
\end{align*}
where $C_{-1}(X):=L_*^{an}(X)$, and $C_{1}(X):=I\chi(X)$.
\end{proof}

\subsection{Proof of Theorem \ref{bsyconjforsymmetric}}
\begin{proof} It suffices to prove the theorem when $X$ is connected.
 Since $X$ is a complex projective variety it is a compact oriented even-dimensional triangulated Witt pseudomanifold so Theorem \ref{symLclassthm} can be applied.
    By the assumption $L_*(X)=\Psi_2IT_{1*}(X)$, since $\Psi_2 L^{an}_*(-)=L_*(-)$ and $IT_{-1*}(X)=I\chi(X)$ by (\ref{mhc-1trivial}), the right hand side of (\ref{symfactorLclass}) equals that of (\ref{symfactorIT1class}). Hence $\Psi_2IT_{1*}(X^{(n)})=\Psi_2 L^{an}_*(X^{(n)})=L_*(X^{(n)})$.
\end{proof}

{\bf AI Declaration.} Prior to submission, the manuscript was reviewed using ChatGPT (GPT-5.6, OpenAI) to identify grammatical  and typographical errors. Some of the references used such as \cite{SV,CEKN,HW,Sag} were suggested to the author by this model. No AI tools were used to develop or verify mathematical arguments or proofs.

\end{document}